\documentclass[11pt, epsfig]{article}
\usepackage{epsfig, amsmath, amssymb, amsthm, times}
\usepackage{setspace}

\renewcommand\Re{\operatorname{Re}}

\newtheorem {thm}{Theorem}[section]
\newtheorem {lem}[thm]{Lemma}

\newtheorem {cor}[thm]{Corollary}

\theoremstyle{defintion}
\newtheorem {df}[thm]{Definition}

\theoremstyle{remark}
\newtheorem{rem}[thm]{Remark}

\theoremstyle{example}
\newtheorem{ex}[thm]{Example}

\theoremstyle{assumption}
\newtheorem{as}[thm]{Assumption}

\def\pf{{\it Proof.\;}}

\def\cov{\mathrm{cov~}}
\def\span{\mathrm{span~}}
\def\tr{\mathrm{Tr~}}
\def\ker{\mathrm{ker~}}
\def\diag{\mathrm{diag~}}

\def\rank{\mathrm{rank~}}

\def\E{{\mathbb E}}

\def\R{{\mathbb R}}
\def\N{{\mathbb N}}
\def\C{{\mathbb C}}

\def\lbl{\label}
\def\be{\begin{equation}}
\def\ee{\end{equation}}
\def\qed{$\square$}
\def\one{\mathbf{1}}
\def\tl{\tilde}
\def\t{\mathsf{T}}

\title{Stochastic stability of continuous time 
consensus protocols}
\author{Georgi S. Medvedev
\thanks{
Department of Mathematics, Drexel University, 3141 Chestnut Street,
Philadelphia, PA 19104, USA, {\tt medvedev@drexel.edu}, 
ph.: 1-215-895-6612, fax: 1-215-895-1582} }
\date{\today}
\begin{document}
\maketitle

\begin{abstract}
A unified approach to studying convergence and stochastic 
stability of  continuous time consensus protocols (CPs) is presented 
in this work. Our method applies to 
networks with directed information flow; both cooperative and
noncooperative interactions; networks under weak stochastic forcing;
and those whose topology and strength of connections may vary in time. 
The graph theoretic interpretation of the analytical results is 
emphasized. We show how the spectral properties, 
such as algebraic connectivity and total effective resistance, as well 
as the geometric properties, such the  dimension and the structure of the 
cycle subspace of the underlying graph, shape
stability of the corresponding CPs.
In addition, we explore certain implications of the spectral 
graph theory to CP design.
In particular, 
we point out that expanders, sparse highly connected graphs, 
generate CPs whose performance remains uniformly high when the size of the
network grows unboundedly. Similarly, we highlight
the benefits of using random versus regular network topologies for 
CP design.
We illustrate these observations with numerical
examples and refer to the relevant graph-theoretic results. 
\end{abstract} 

\section{Introduction}
\setcounter{equation}{0}
The theory of CPs is a framework for design and analysis of distributed 
algorithms for coordination of the groups of dynamic agents. 
In many control problems, agents in the group need to agree upon certain
quantity, whose interpretation depends on the problem at hand. The theory
of CPs studies the convergence to a common value (consensus) in its general
and, therefore, abstract form.
It has been
a subject of intense research due to diverse applications in applied science
and engineering. The latter include coordination of groups of unmanned 
vehicles \cite{Ren07, Olfati07}; synchronization of power, sensor and 
communication networks \cite{CB05,DB10}; and principles underlying collective 
behavior in social networks \cite{Olfati07a} and biological systems
\cite{balerini, sumpter}, to name a few. 

From the mathematical point of view, analysis of continuous time CPs
is a stability problem for systems of linear differential equations
possibly with additional features such as stochastic perturbations
or time delays. There are many effective techniques for studying
stability of linear systems \cite{Demidovich, hale_odes}.
The challenge of applying
these methods to the analysis of CPs is twofold. First, one is interested
in characterizing stability under a minimal number of practically relevant
assumptions on the structure of the matrix of coefficients, which may
depend on time. Second, it is important to identify the relation
between the structure of the graph of interactions
in the network to the dynamical performance of CPs. 
A successful solution of the second problem requires
a compilation of dynamical systems and graph theoretic
techniques. This naturally leads to spectral methods, 
which play important roles in both mathematical disciplines, 
and are especially useful for problems on the interface between 
dynamics and the graph theory \cite{Jost07}. A general idea for using 
spectral methods for analyzing CPs is that, on the one hand, 
stability of the continuous time CP is encoded in the eigenvalues (EVs) 
of the matrix of coefficients; on the other hand, EVs of the same
matrix capture structural properties of the graph of the CP.
The spectral graph theory offers many fine results relating 
the structural properties of graphs to 
the EVs of the adjacency matrix and the graph Laplacian
\cite{Fiedler73, Allon86, Chung97, Cioaba10, Hoory06, Nil91, Gutman03}. 
This provides a link between the network topology and the dynamical
properties of CPs.

In this paper, under fairly general assumptions on CPs, we study two
problems: convergence of CPs and their stability in the presence of
stochastic perturbations. The former is the problem of asymptotic stability
of the consensus subspace, a one-dimensional invariant (center) subspace.
The latter is a special albeit representative form of stability of the
consensus subspace with respect to constantly acting perturbations
\cite{Malkin}.
The rate of convergence to the consensus
subspace sets the timescale of the consensus formation (or synchronization)
from arbitrary initial conditions or upon instantaneous perturbation.
Therefore, the convergence rate is important in applications where
the timing of the system's responses matters (e.g., in decision making 
algorithms, neuronal networks, etc). Stochastic stability, on the other
hand, characterizes robustness of the consensus to noise. This form of stability
is important when the consensus needs to be maintained in noisy environment
over large periods of time (e.g., communication networks, control of
unmanned vehicles, etc). We believe that our quantitative description of these
two forms of stability elucidates two important aspects of the performance
of CPs.

The questions investigated in this paper have been studied before under
various hypotheses on CPs: constant weights \cite{Olfati04, Olfati07, Ren07}, 
time-dependent interactions \cite{Mor04, RB05, Ren07}, and CPs with time-delays 
\cite{Lin08, Olfati04}. Optimization problems arising in the context of 
CP design were studied in \cite{Boyd06, Boyd03, Boyd06a}. 
There is a body of related work on discrete time CPs \cite{Boyd03, Olfati07}. 
Robustness of CPs  to noise was studied in \cite{Leonard10, LZ09}.
In this paper, we offer a unified approach to studying convergence and
stochastic stability of CPs. 
Our method applies to networks with directed information flow; 
both cooperative and noncooperative interactions; networks under 
weak stochastic forcing; and those whose topology and strength of 
connections may vary in time. We derive sufficient conditions 
guaranteeing convergence of time-dependent CPs and present estimates
characterizing their stochastic stability. For CPs on undirected 
graphs, we show that the rate of convergence and stability
to random perturbations are captured by the generalized algebraic connectivity
and the total effective resistance of the underlying graphs. 
Previously, these results were available only for CPs on graphs
with positive weights  \cite{Mor04, Leonard10}.
To further elucidate the role that network topology plays in shaping
the dynamical properties of CPs, we further develop our results for CPs
on simple networks (see text for the definition of a simple network).
Our analysis of simple networks reveals the role of the geometric properties 
of the cycle subspace associated with the graph of the network (such as
the first Betti number of the graph; the length and the mutual position
of the independent cycles) to the stability of CPs to random perturbations.
In addition, we explore several implications of the results of the spectral 
graph theory to CP design. First, we show that expanders, sparse highly connected
graphs \cite{Hoory06,Sar04}, generate CPs with the rate of convergence
bounded from zero uniformly when the size of the network tends to infinity.
In particular, CPs based on expanders are effective for coordinating
large networks. Second, we point out that CPs with random connections have
nearly optimal convergence rate. In contrast, the convergence of CPs on regular 
lattice-like graphs slows down rapidly as the size of
the network grows. We illustrate these observations with numerical
examples and refer to the relevant graph-theoretic results. 

The mathematical analysis of CPs in this paper uses the method, 
which we recently developed for studying 
synchronization in systems of coupled nonlinear oscillators 
\cite{medvedev10, medvedev10a} and reliability of neuronal networks
\cite{medvedev09}. We further develop this method in several ways. First,
we relate the key properties of the algebraic 
transformation of the coupling operator used in 
\cite{medvedev10, medvedev10a, medvedev09} for studying synchronization
to general properties of a certain class
of pseudo-similarity transformations. Second, we strengthen the graph 
theoretic interpretation of the stability analysis. 
We believe  that our method will be useful for design and analysis of CPs and 
for studying synchronization in a large class of models.

The outline of the paper is as follows. In Section~\ref{algebra},
we study the properties of a pseudo-similarity transformation, which
is used in the analysis of CPs in the remainder of the paper.
Section~\ref{convergence} is devoted to the convergence analysis
of CPs. After formulating the problem and introducing 
necessary terminology in \S\ref{formulation}, we study convergence
of CPs with constant and time-dependent coefficients in \S\ref{stationary}
and \S\ref{nonstationary} respectively. Section~\ref{robustness} presents
estimates characterizing stochastic stability of stationary and time-dependent
CPs. These results are applied to study CPs protocols on undirected weighted
graph in Section~\ref{undirected}.
In Section~\ref{connectivity}, we discuss the relation between the
connectivity of the graph and dynamical performance of CPs. 
The results of this paper are summarized in Section~\ref{conclude}.

\section{Pseudo-similarity transformation}\lbl{algebra}
The analysis of CPs in the sections that follow relies on certain properties
of a pseudo-similarity transformation, which we study first.  
\begin{df}\lbl{pseudo-sim} 
Matrix $\hat D\in \R^{(n-p)\times (n-p)} \; (1\le p\le n-1)$
is pseudo-similar to $D\in\R^{n\times n}$ via $S\in\R^{(n-p)\times n}$ if
\be\lbl{com}
SD=\hat DS.
\ee
\end{df}
Equation (\ref{com}) is equivalent to the following property
\be\lbl{poly}
Sq(D)=q(\hat D)S,
\ee
for any polynomial $q(t)$.

To study the existence and the properties of pseudo-similar matrices,
we recall the definition of the Moore-Penrose pseudo-inverse 
of a rectangular matrix (cf. \cite{Gan}).
\begin{df}\lbl{Moore}
$A^+\in\R^{n\times m}$ is called a pseudo-inverse of $A\in\R^{m\times n}$ if 
$$
AA^+A=A\quad\mbox{and}\quad A^+AA^+=A^+.
$$
\end{df}  

Throughout this section, we use the following assumption.
\begin{as}\lbl{defineS}
Let $D\in\R^{n\times n}$ and $S\in\R^{(n-p)\times n}\; (1\le p\le n-1)$ such that
\begin{eqnarray}
\lbl{rankS}
\rank S &=& n-p,\\
\lbl{kerS}
\ker S &\subset& \ker D.
\end{eqnarray}
\end{as}
Condition (\ref{rankS}) implies that 
\be\lbl{full-row}
S^+=S^\t  (SS^\t )^{-1},
\ee
and, therefore,
\be\lbl{two-prop}
S^+S=P_{R(S^\t )}=P_{(\ker S)^\perp}, \quad\mbox{and}\quad SS^+=I_{n-p}.
\ee
Here, $P_{R(S^\t )}$ and $I_{n-p}$ denote the projection matrix
onto the column space of $S^\t $ and the $(n-p)\times (n-p)$ identity
matrix. 

The combination of (\ref{rankS}) and (\ref{kerS}) guarantees the
existence and uniqueness of the pseudo-similar matrix for $D$
via $S$. 
\begin{lem}\lbl{main-property}
Let $S\in\R^{(n-p)\times n}$ and $D\in\R^{n\times n}$ satisfy Assumption~\ref{defineS}. 
Then 
\be\lbl{pseudo}
\hat D=SDS^+
\ee
is a unique pseudo-similar matrix to $D$ via $S$.
\end{lem}
\pf 
By the first identity in (\ref{two-prop}),
$$
DS^+S=D.
$$
Therefore, equation (\ref{com}) is solvable with respect to 
$\hat D$. By multiplying both sides of (\ref{com}) by
$S^+$ from the right and using the second property in (\ref{two-prop}),
we obtain (\ref{pseudo}).\\
\qed
\begin{cor}\lbl{two-more}
\be\lbl{exp}
 \exp\{t\hat D\}= S\exp\{tD\}S^+,\quad t\in\R.
\ee
\end{cor}
\pf 
Equation (\ref{exp}) follows from the second identity in (\ref{two-prop})
and the series representation of $\exp\{tD\}$.\\
\qed

The next lemma relates the spectral properties of $D$ and $\hat D$.
\begin{lem}\lbl{spectra}
Suppose $D$ and $S$ satisfy Assumption~\ref{defineS} and $\hat D$ is
the pseudo-similar matrix to $D$ via $S$.
\begin{description}
\item[A] If $\lambda\in\C$ is a nonzero EV of $D$
then $\lambda$ is an EV of $\hat D$ of the same algebraic
and geometric multiplicity. Moreover, $S$ carries out 
a bijection from the generalized $\lambda-$eigenspace
of $D$ onto that of $\hat D$ preserving the Jordan block
structure.
\item[B]
$\lambda=0$ is an EV of $\hat D$ if and only if the algebraic multiplicity of
$0$ as an EV of $D$ exceeds $p$. In this case, the algebraic
multiplicity of $0$ as  as an EV of $\hat D$ is diminished by
$p$. $S$ maps the generalized $0-$eigenspace of $D$ onto that
of $\hat D$.
\item[C] $S$ maps a Jordan basis of $D$ onto that of $\hat D$.
\end{description}
\end{lem}
\pf 
\begin{description}
\item[A]
$S$ restricted to the direct sum of generalized eigenspaces of $D$ 
corresponding to nonzero eigenvalues is injective.

Let $\lambda$ be a nonzero EV of $D$.  Since for any $m\in\N,$
$
(\hat D-\lambda I_{n-p})^mS =S(D-\lambda I_n)^m
$
(cf. (\ref{poly})), $Sv$ is a generalized $\lambda-$eigenvector of $\hat D$ of index $m$
if and only if $v$ is a generalized $\lambda-$eigenvector of $D$ of index $m$.
Therefore, $S$ bijectively maps the generalized $\lambda-$eigenspace of $D$
onto that of $\hat D$. The associated Jordan block structures are the
same.
\item[B]
If the generalized $0-$eigenspace of $D$ is larger than $\ker S$ 
then $\ker\hat D$ is nontrivial.
Choose a Jordan basis for $D$ restricted to its generalized $0-$eigenspace
\begin{eqnarray*}
v_1^{(1)}, v_2^{(1)},\dots v_{k_1}^{(1)}, & (\ker D^m\ominus \ker D^{m-1})\\
Dv_1^{(1)}, Dv_2^{(1)},\dots Dv_{k_1}^{(1)},v_1^{(2)}, v_2^{(2)},\dots v_{k_2}^{(2)},
& (\ker D^{m-1}\ominus \ker D^{m-2})\\
\dots & \\
D^{(m-1)}v_1^{(1)}, D^{(m-1)}v_2^{(1)},\dots D^{(m-1)}v_{k_1}^{(1)},
\dots, v_1^{(m)}, v_2^{(m)},\dots v_{k_m}^{(m)}
& (\ker D).
\end{eqnarray*}
The image of this basis under $S$ consists of the vectors forming
a Jordan basis of $\hat D$ restricted to its generalized $0-$eigenspace
and $p$ zero vectors. Under the action of $S$, each cyclic subspace
of $D$
$$
\span\left(v_j^{(i)}, Dv_j^{(i)}, \dots D^{m-i}v_j^{(i)}\right)
$$
looses a unit in dimension if and only if  $D^{m-i}v_j^{(i)}\in \ker S$.

\item[C]
The statement in \textbf{C} follows by applying the
argument in \textbf{B} to a Jordan basis of $D$ restricted to 
the generalized eigenspace corresponding to a nonzero eigenvalue.
\end{description}
\qed

Next, we apply Lemmas~\ref{main-property} and \ref{spectra} to 
the situation, which will be used in the analysis of CPs below.
\begin{cor}\lbl{1dker}
Denote $e=(1,1,\dots,1)^\t\in\R^n$ and $\one=\span\{e\}$. 
Let $D\in\R^{(n-1)\times n}$ and $S\in\R^{(n-1)\times n}$ be such that
\be\lbl{D-and-S}
D\in\mathcal{K}=\{M\in\R^{n\times n}:~ Me=0\}\;\mbox{and}\;
\ker S=\one.
\ee

By Lemmas~\ref{main-property} and \ref{spectra}, we have
\begin{enumerate}
\item
\be\lbl{exists}
\exists !~\hat D\in\R^{(n-1)\times (n-1)}:~ SD=\hat DS
\ee
\item
\be\lbl{Dhat}
\hat D=SDS^+.
\ee
\item
Denote the EVs of $D$ counting multiplicity by 
\be\lbl{specD}
\lambda_1=0, \lambda_2, \lambda_3, \dots, \lambda_{n},
\ee
such that $e$ is an eigenvector corresponding to $\lambda_1$.
Then  
\be\lbl{specDhat}
\lambda_2, \lambda_3, \dots, \lambda_{n}
\ee
are the EVs of $\hat D$. For $i=2,3,\dots,n,$ $S$ maps bijectively 
the generalized $\lambda_i-$eigenspaces of $D$ to those of $\hat D$.
$0$  is an EV of $\hat D$ if and only if the algebraic multiplicity
of $0$ as an EV of $D$ is greater than $1$.
\end{enumerate}
\end{cor}
\begin{ex}\lbl{main-example} The following matrix
satisfies (\ref{D-and-S}) and can be used as an intertwining 
matrix in (\ref{exists})
\be\lbl{matrixSexample}
S=
\left(\begin{array}{cccccc}
-1 & 1 & 0& \dots &0& 0 \\
0 & -1& 1 & \dots& 0& 0\\
\dots&\dots&\dots&\dots&\dots&\dots\\
0 &0 & 0& \dots &-1 & 1
\end{array}
\right) \in\R^{(n-1)\times n}.
\ee
\end{ex}
\section{Convergence analysis of CPs}\lbl{convergence}
\setcounter{equation}{0}
In \S\ref{formulation}, we introduce a continuous time CP, a differential 
equation model that will be studied in the remainder of this paper.
Convergence of CPs with constant and 
time-dependent coefficients is analyzed in \S\ref{stationary} and
\S\ref{nonstationary}, respectively.

\subsection{The formulation of the problem}\lbl{formulation}

By a continuous time CP with constant coefficients we call
the following  system of ordinary differential equations (ODEs):
\be\lbl{CP}  
\dot x^{(i)}=\sum_{j=1}^n a_{ij} (x^{(j)}-x^{(i)}),\;\; i\in[n]:=\{1,2,\dots,n\}.
\ee
Unknown functions $x^{(i)}(t), i\in [n]$ are interpreted as
the states of $n$ agents. The right hand
side of (\ref{CP}) models the information exchange between agents
in the network. 
Coefficient $a_{ij}$ is interpreted as the weight
that Agent~$i$ attributes to the information from Agent~$j$. 
Positive weights promote synchronization between the states of 
the corresponding agents. Negative weights have the opposite 
effect and can be used
to model noncooperative interactions between the agents. For more 
background and motivation for considering (\ref{CP}) and related 
models we refer an interested reader to 
\cite{Mor04, Olfati07a, Olfati04, Ren07, Leonard10}.

An important problem in the analysis of CPs is identifying conditions,
under which the states of the agents in the network converge to 
the same value.
\begin{df} We say that CP (\ref{CP}) reaches a consensus from 
initial state $x(0)\in\R^n$ if 
\be\lbl{consensus}
\lim_{t\to\infty} |x^{(j)}(t)-x^{(i)}(t)|=0 \;\forall (i,j)\in [n]^2.
\ee
If CP (\ref{CP}) reaches a consensus from any initial 
condition then it is called convergent.
\end{df}
The second problem, considered in this paper, is that of stability of the consensus subspace $\one$
to instantaneous and constantly acting random perturbations.
An important aspect of the analysis of convergence and stability of CPs 
is elucidating the relation between the structural properties of the 
network (e.g., connectivity and weight distribution) and the degree of
stability of the corresponding CPs. To study this problem, we need
the following graph-theoretic interpretation of (\ref{CP}).

Let $A=(a_{ij})\in\R^{n\times n}$ be a matrix of coefficients of (\ref{CP}).
Since (\ref{CP}) is independent from the choice of diagonal elements
of $A$, we set $a_{ii}=0, i\in [n].$ Using the terminology from the electrical networks
theory \cite{Bollobas98}, we call $A$ a conductance matrix. Next, we associate
with (\ref{CP}) a directed graph $G=(V,E)$, where the vertex set $V=[n]$
lists all agents and a directed edge $(i,j)\in [n]^2$ belongs to $E$ if 
$a_{ij}\ne 0.$ 

By the  network we call 
$
\mathcal{N}=(V, E, a)= (G,a),
$
where function $a: E\to \R$ assigns conductance $a_{ij}\in\R$ to each edge
$(i,j)\in E$. If conductance matrix $A$ is symmetric, $G$ can be viewed
as an undirected graph. If, in addition, $a_{ij}\in\{0,1\}$, $\mathcal{N}$
is called simple.


\subsection{Stationary CPs}\lbl{stationary}
The convergence analysis of CPs with constant and time-dependent coefficients
relies on standard results of the theory of differential equations
(see, e.g., \cite{hale_odes}). It is included for completeness and to 
introduce the method that will be used later for studying stochastic 
stability of CPs.

We rewrite (\ref{CP}) in matrix form
\be\lbl{constant}
\dot x = Dx.
\ee
The matrix of coefficients
\be\lbl{coupling}
D=A-\diag (\bar a_{1},\bar a_{2},\dots,\bar a_{N}),\quad
\bar a_{i}=\sum_{j=1}^N a_{ij}
\ee
is called a coupling matrix.

Let $\tilde S\in\R^{(n-1)\times n}$ be a matrix with one
dimensional null space, $\ker\tilde S=\one$  
(see Example~\ref{main-example} for a possible choice of $\tl S$).  The analysis
in this section does not depend on the choice of $\tl S$.
Suppose $\tl S$ has been fixed and define
\be\lbl{rescaleS}
S=(\tl S{\tl S}^\t)^{-{1\over 2}}\tl S.
\ee
Note that $S$ has orthogonal rows
\be\lbl{orthogonality}
SS^\t=I_{n-1} \quad\mbox{and}\quad S^\t S= \tl S^+ \tl S=P_{\one^\perp},
\ee
where $P_{\one^\perp}$ stands for the orthogonal projection onto
$\one^\perp$.
By definition, $D$ and $S$ satisfy conditions of Corollary~\ref{1dker}.
Therefore, there exists a unique $S-$reduced matrix
\be\lbl{newDhat}
\hat D=SDS^+=(\tl S{\tl S}^\t)^{-{1\over 2}}\tl SD{\tl S}^\t
(\tl S{\tl S}^\t)^{-{1\over 2}}=SDS^\t,
\ee
whose properties are listed in Corollary~\ref{1dker}.
In addition, using normalized matrix $S$ (cf. (\ref{rescaleS}))
affords the following property. 

\begin{lem}\lbl{normal}
Let $D\in\mathcal{K}$ and $S$ be as defined in (\ref{rescaleS}).
Then $\hat D$, the pseudo-similar matrix to $D$ via $S$, is normal (symmetric) 
if $D$ is normal (symmetric).
\end{lem}
\pf
If $D$ is symmetric, then so is $\hat D$ by (\ref{newDhat}).

Suppose $D$ is normal. Then there exists an orthogonal matrix
$
V=(e, v_2,v_3,\dots, v_n)
$ 
and diagonal matrix 
$$
\Lambda=\mathrm{diag}\left(0,\lambda_2,\dots,\lambda_{k_1},
\begin{pmatrix} 
\alpha_{k_1+1} & -\beta_{k_1+1}\\\beta_{k_1+1} &\alpha_{k_1+1}
\end{pmatrix},\dots,
\begin{pmatrix} 
\alpha_{k_l}&  -\beta_{k_l}\\ \beta_{k_l} &\alpha_{k_l}
\end{pmatrix}\right).
$$
such that 
$
D=V\Lambda V^\t.
$
By Lemma~\ref{spectra}, 
$
D=U\hat\Lambda U^\t.
$
with 
$
U=(Sv_2,Sv_3,\dots, Sv_n)
$ 
and 
$$
\hat\Lambda=\mathrm{diag}\left(\lambda_2,\dots,\lambda_{k_1},
\begin{pmatrix} 
\alpha_{k_1+1} & -\beta_{k_1+1}\\\beta_{k_1+1} &\alpha_{k_1+1}
\end{pmatrix},\dots,
\begin{pmatrix} 
\alpha_{k_l}&  -\beta_{k_l}\\ \beta_{k_l} &\alpha_{k_l}
\end{pmatrix}\right).
$$
Denote the columns of $U$ by $u_i$, $i\in [n]/\{1\}$.
Since
$$
u^\t_j u_i= v_j^\t S^\t S v_i= v_jP_{1^\perp}v_i =v_j^\t v_i, 
\; i,j \in\{2,3,\dots,n\},
$$
$U$ is an orthogonal matrix. Therefore, $\hat D$ is normal.\\
\qed

By multiplying both sides of (\ref{constant}) by
$S$ and using (\ref{com}), we obtain the reduced equation
for $y=Sx\in R^{n-1}:$
\be\lbl{reduced}
\dot y =\hat Dy.
\ee
Under $S$, the consensus subspace of (\ref{constant}) is mapped to
the origin of the phase space of the reduced system. Furthermore,
because $\hat D$ inherits its spectrum from $D$ (cf. Lemma~\ref{spectra}), there is 
a direct relation between the transverse stability of the consensus subspace
$\one$, as an invariant center subspace of the original problem 
(\ref{CP}), and that of the equilibrium at the origin of the reduced
system (\ref{reduced}). This relation is described in the following theorem.

\begin{thm}\lbl{exhaustive}
CP (\ref{constant}) is convergent if and only if 
\be\lbl{stable}
D\in\mathcal{D}:=\{M\in\mathcal{K}\;\&\;\hat M\; 
\mbox{is a stable matrix}\},
\ee
where $\hat{M}$ is the pseudo-similar matrix to $M$ via $S$ 
(cf. (\ref{rescaleS})).
If $D\in\mathcal{D}$, the rate of convergence to the consensus
subspace is set by the nonzero  EV of $D$ with the largest
real part. Specifically, let the EVs of $D$ be arranged as
in Corollary~\ref{1dker} and
\be\lbl{alpha} 
\alpha =-\max_{i\ge 2} \Re\lambda_i.
\ee
Then there exists $C_1>0$ such that for any initial
condition $x(0)\in\R^n$ and any $\epsilon>0$
\be\lbl{rate}
|P_{\one^\perp}x(t)|\le C_1 |P_{\one^\perp}x(0)|\exp\{(-\alpha+\epsilon) t\}, 
\ee
where $|\cdot|$ stands for the Euclidean norm in $R^{n-1}$. 
\end{thm}
\pf
Let $x(t)$ be a solution of (\ref{constant}). Denote the projection
of $x(t)$ onto $\one^\perp$ by 
$$
z(t)=P_{\one^\perp}x(t)=S^\t S x(t).
$$
By multiplying both parts of (\ref{constant}) by $S^\t S$ and using
(\ref{orthogonality}), we derive an ode for $z(t)$
\be\lbl{eqn-for-z}
\dot z= S^\t \hat DSz.
\ee
On the other hand, 
$$
y(t)=Sx(t)=SS^\t Sx(t)=Sz(t)
$$
satisfies the reduced equation $\dot y= \hat D y$.
Therefore, $S:\R^n\to\R^{n-1}\cong \one^\perp$ provides a one-to-one
correspondence between the trajectories of the reduced system and the 
projections of the trajectories of (\ref{constant}) onto the orthogonal
complement of the consensus subspace, $\one^\perp$. In addition, $S$ maps
the consensus subspace to the origin of phase space of the reduced system.
Therefore, transverse asymptotic
stability of $\one$ as an invariant linear subspace of (\ref{constant})
is equivalent to the asymptotic stability of the fixed point at the 
origin of the reduced system. The necessary and sufficient condition for the latter
is that $\hat D$ is a stable matrix, i.e., $D\in\mathcal{D}$.

If $D\in\mathcal{D}$ then $\lambda_1=0$ is a simple eigenvalue
and the real parts of the remaining EVs are negative. 
By standard results of the theory of ODEs, there exists a positive 
constant $C_1$ such that for any initial
condition $y(0)=Sx(0)\in\R^{n-1}$ and any $\epsilon>0$, the solution
of the reduced system satisfies
\be\lbl{rate-for-y}
|y(t)|\le C_1 |y(0)|\exp\{(-\alpha+\epsilon) t\}. 
\ee
Therefore,
$$
|P_{\one^\perp} x(t)=
|S^\t y(t)|=|y(t)|\le C_1 |y(0)|\exp\{(-\alpha+\epsilon) t\}. 
$$
\qed

For CPs with nonnegative weights $a_{ij}\ge 0$ there is a simple
sufficient condition of convergence: nonnegative CP (\ref{constant})
is convergent if the corresponding digraph is strongly connected
\cite{Olfati04}. This condition does not hold in general if there are 
negative weights. Edges with positive weights help
synchronization. In contrast, negative weight $a_{ij}<0$ indicates that 
Agent~$i$ does not cooperate with Agent~$j$ in reaching consensus.
The lack of cooperation by some agents in the network can be compensated 
by the cooperative interactions between other agents. The consensus
can be reached even if many of the weights are negative, as shown
in the following example. 

\begin{ex}~\lbl{uniform}
The following example of a random matrix $D\in\mathcal{D}$
was constructed in \cite{medvedev10}:
\be\lbl{random_D}
D=
\left(\begin{array}{ccccc}
 -1.0251   & 2.2043 &  -1.6032 &   0.5044   & -0.0804\\
   -0.1264 &   0.2772 &  -0.3006 &   0.2060 &   -0.0562\\
   -1.1549 &   2.5819 &  -1.9613 &   0.5210 &   0.0133 \\
   -0.8807 &   1.9231 &  -1.0823 &   0.0333 &   0.0066 \\
   -0.9049 &   1.8778 &  -1.0060 &   0.3772 &  -0.3441
\end{array}
\right).
\ee
About half of the entries of $D$ 
are negative, i.e., there are as many noncooperative interactions
as cooperative ones. Nonetheless, the resultant CP is convergent.
\end{ex}

\subsection{Time-dependent CPs}\lbl{nonstationary}
In realistic networks, the strength and even topology 
of connections between the agents may depend on time. To account for
a greater variety of possible modeling situations, in this section,
we use only very mild assumptions on the regularity of conductances 
$a_{ij}(t)$ as functions of time: 
$a_{ij}(t), (i,j)\in [n]^2,$ are measurable 
locally bounded real-valued functions.
Under these assumptions, we formulate
two general sufficient conditions for convergence of time-dependent
CPs. 

By a time-dependent CP, we call the following ODE
\be\lbl{dynamic}
\dot x = D(t)x,\; D(t)=(d_{ij}(t))\in\R^{n\times n},
\ee
where the coupling matrix $D(t)$ is defined as before
\be\lbl{coupling-again}
D(t)=A(t)-\diag(\bar a_1(t),\bar a_2(t),\dots, \bar a_n(t)), \;
\bar a_i(t)=\sum_{j=1}^n a_{ij}(t).
\ee 
Under our assumptions on $a_{ij}(t)$, the solutions of (\ref{dynamic})
(interpreted as solutions of the corresponding integral equation) are
well-defined (cf. \cite{hale_odes}).

By construction, coupling matrix $D$ satisfies the following condition
(see (\ref{coupling}))
$$
D(t)\in \mathcal{K}\;\forall t\ge 0.
$$
For convenience, we reserve a special notation for the class
of admissible matrices of coefficients:
\be\lbl{admissible}
D(t)\in\mathcal{K}_1:=\{
M(t)\in\R^{n\times n}:\; m_{ij}(t)\in L_{loc}^\infty(\R^+)\;\&\;
(M(t)\in\mathcal{K}\;\forall t\ge 0)
\}.
\ee
Note that for any $t>0$ there exists a unique pseudo-similar matrix to
$D(t)$ via $\tilde S$, $\hat D(t)=\tilde SD(t){\tilde S}^+$, provided 
$\tilde S$ satisfies Assumption~\ref{defineS}.
Below, we present two classes of 
convergent time-dependent CPs.
The first class is motivated by the convergence analysis of CPs with
constant coefficients.
\begin{df}\lbl{dissipative} \cite{medvedev10}
Matrix valued function $D(t)\in\mathcal{K}_1$ is called uniformly 
dissipative with parameter $\alpha$ if there exists $\alpha>0$ such that 
\be\lbl{uniform}
y^\t  \hat D(t) y \le -\alpha y^\t  y \quad\forall y\in\R^{n-1}\;\&\; \forall t\ge 0,
\ee
where $\hat D(t)$ is the pseudo-similar matrix to $D(t)$ via $\tilde S$.
The class of uniformly dissipative matrices is denoted by $\mathcal{D}_\alpha$.
\end{df}
The convergence of uniformly dissipative CPs is given in the following theorem.
\begin{thm}\lbl{uniform_convergence}
Let $D(t)\in\mathcal{D}_\alpha,\;\alpha>0$.
CP (\ref{dynamic}) is convergent with the rate of convergence
at least $\alpha$.
\end{thm}
\pf
It is sufficient to show that $y(t)\equiv 0$ is an asymptotically stable
solution of the reduced system for $y=\tilde Sx$
\be\lbl{reduced-again}
\dot y=\hat D(t).
\ee
For solution $y(t)$ of the reduced system, we have
$$
{d\over dt} |y|^2=2 y^\t \hat D(t)y\le 2 \alpha |y|^2.
$$
Thus,
$$
|y(t)|\le |y(0)| \exp\{\alpha t\}.
$$
\qed

In conclusion, we prove convergence for a more general class of CPs.
\begin{df}\lbl{asympt-dissipative}
The coupling matrix $D(t)$ is called asymptotically dissipative if
\be\lbl{average}
D\in\mathcal{\tl D}=\{ M(t)\in\mathcal{K}_1:\; 
\limsup_{t\to\infty} t^{-1}\int_0^t \sup_{|y|=1} y^\t  \hat M(u) y du <0 \}.
\ee
\end{df}
\begin{thm}\lbl{Vazh}
If $D(t)\in\mathcal{\tl D}$ then CP (\ref{dynamic}) is convergent.
\end{thm}
\pf 
Let $y(t)$ be a solution of the reduced system (\ref{reduced-again}).
Then
$$
{d\over dt} |y|^2=2 y^\t  \hat D(t)y\le 2 \gamma(t) |y|^2,\;\;
\gamma(t):=\sup_{|y|=1} y^\t \hat D(t)y.
$$
By Gronwall's inequality,
$$
|y(t)|\le |y(0)| \exp\{\int_0^t\gamma(u)du\}.
$$
Since $D(t)\in \mathcal{\tl D}$, 
$$
\exists \alpha>0\; \&\; T>0\;:\; (t\ge T)\Rightarrow \int_0^t\gamma(u)du\le -\alpha t.
$$
Thus,
$$ 
|y(t)|\le |y(0)| \exp\{-\alpha t\}< -\alpha t,\; t\ge T.
$$
\qed

\section{Stochastic stability}\lbl{robustness}
\setcounter{equation}{0}
In this section, we study stochastic stability of CPs.
Specifically, we  consider
\be\lbl{perturbed}
\dot x = D(t)x+\sigma U(t)\dot w,\;\; x(t)\in\R^n,
\ee
where $\dot w$ is a white noise process in $\R^n$, 
$D(t)\in\mathcal{K}_1$ (cf. (\ref{admissible})), and $U(t)=(u_{ij}(t))\in\R^{n\times n}$,
$u_{ij}(t)\in L^\infty_{loc}(\R^+).$
The consensus subspace, $\one$, forms an invariant center subspace
of the deterministic system obtained from (\ref{perturbed}) by setting $\sigma=0$.

Since the transverse stability of the consensus subspace is equivalent 
to the stability of the equilibrium of the corresponding 
reduced equation, along with
(\ref{perturbed}), we consider the corresponding equation for $y=Sx$:
\be\lbl{red}
\dot y=\hat D(t)y+\sigma SU(t)\dot w,
\ee 
where $S$ is defined in (\ref{rescaleS}).
The solution of (\ref{red}) with deterministic initial condition 
$y(0)=y_0\in\R^{n-1}$
is a Gaussian random process.
The mean vector and the covariance matrix functions of stochastic process $y(t)$
\be\lbl{rho-and-e}
m(t):= \E y(t) \quad \mbox{and}\quad 
V(t):=\E \left[(y(t)-m(t))(y(t)-m(t))^\t\right],
\ee
satisfy linear equations (cf. \cite{KS})
\be\lbl{linear}
\dot m= \hat D m\quad \mbox{and}\quad 
\dot V= \hat D V +V{\hat D}^\t+\sigma^2 SU(t)U(t)^\t S^\t.
\ee
The trivial solution  $y(t)\equiv 0$ of the reduced equation 
(\ref{reduced}) is not a solution of the perturbed equation
(\ref{red}).  Nonetheless,
if the origin is an asymptotically stable equilibrium
of the deterministic reduced equation obtained from (\ref{red})
by setting $\sigma=0$, the trajectories of (\ref{perturbed})
exhibit stable behavior. In particular, if (\ref{dynamic}) is
a convergent CP, for small $\sigma>0$,
the trajectories of the perturbed system (\ref{perturbed})
remain in $O(\sigma)$ vicinity of the consensus subspace
on finite time intervals of time with high probability. We use
the following form of stability to describe this situation formally.
\begin{df}\lbl{stability}
CP (\ref{dynamic}) is stable to
random perturbations (cf. (\ref{perturbed})) if for any initial 
condition $x(0)\in\R^n$ and $T>0$
\be\lbl{two-moments}
\lim_{t\to\infty} \E P_{\one^\perp}x(t)=0\quad\mbox{and}\quad 
\E|P_{\one^\perp}x(t)|^2=O(\sigma^2),\; 
t\in [0,T].
\ee  
\end{df}
\begin{thm}\lbl{robust}
Let $D\in\mathcal{D}_\alpha$ be a uniformly dissipative matrix
with parameter $\alpha$ (cf. Definition~\ref{dissipative}). 
Then CP (\ref{dynamic}) is stable to random perturbations.
In particular, the solution of the initial 
value problem for (\ref{perturbed}) with deterministic
initial condition $x(0)\in\R^{n}$ satisfies
the following estimate
\be\lbl{main}
\E|P_{\one^\perp}x(t)|^2 
 \le {\sigma^2 n \over 2 \alpha}
\sup_{u\in [0,t]}\Vert U(u)U^\t (u)\Vert,\; t>0,
\ee
where $\Vert \cdot \Vert$ stands for the operator matrix norm 
induced by the Euclidean norm.
\end{thm}
\begin{rem}
Suppose the strength of interactions between the agents in the network
can be controlled by an additional parameter $g$ 
\be\lbl{strength}
\dot x= g D(t)x +\sigma U(t)\dot w.
\ee
Here, the larger values of $g>0$ correspond to stronger coupling, i.e., 
to  faster information exchange in the network. 
By applying estimate (\ref{main}) to
(\ref{strength}), we have
\be\lbl{new_kappa}
\E|P_{\one^\perp}x(t)|^2 \le {\sigma^2 n \over 2 g\alpha}
\sup_{u\in [0,t]}\Vert U(u)U^\t (u)\Vert,\; t>0.
\ee
Note that the variance of $|P_{\one^\perp}x(t)|$ can be effectively 
controlled by $g$. 
In particular, the accuracy of the consensus can be enhanced to any 
desired degree by increasing the rate of information exchange between 
the agents in the network. For the applications of this observations to 
neuronal networks, we refer the reader to \cite{medvedev09, MZ11}. 
\end{rem}
\begin{rem}\lbl{note}
Since $P_{\one^\perp} x(t)=S^\t y(t)$ and  
$(P_{\one^\perp}x(t))^\t P_{\one^\perp} x(t)=y(t)^\t y(t)$,
it sufficient to prove (\ref{main}) 
with $P_{\one^\perp} x(t)$ replaced by the solution of the reduced
problem (\ref{red}), $y(t)=Sx(t)$.
\end{rem}
\pf Let $\Phi(t)$ denote the principal matrix solution of the homogeneous
equation (\ref{red}) with $\sigma=0$.  The solution of the initial value problem 
for (\ref{red}) is a Gaussian random process whose expected value and 
covariance matrix are given by
\begin{eqnarray}\lbl{Ey}
\E y_t &=& \Phi(t)y_0,\\
\lbl{covy}
\cov y_t &=& \sigma^2 
\Phi(t) \int_0^t \Phi^{-1}(u) SU(u)U(u)^\t S^\t  (\Phi^{-1}(u))^\t du~\Phi(t)^\t . 
\end{eqnarray}
Since $D(t)\in\mathcal{D}_\alpha$,
we have
\be\lbl{unid}
y^\t \hat D(t)y\le -\alpha y^\t y \quad\forall y\in\R^{n-1},\; t\ge 0.
\ee
This has the following implication. For all $t\ge 0$,
$-\hat D^s(t)$ is positive definite and, therefore,
$(-\hat D^s(t))^{-1\over 2}$ is well defined. 
Here and throughout this paper, $M^s:= 2^{-1}(M+M^\t)$ stands for
the symmetric part of square matrix $M$.

Using this observation, we rewrite the integrand in (\ref{covy}) as follows
\footnote[1]{The derivations of the estimates (\ref{integrand}) and (\ref{there})
use the following matrix 
inequality: $A(B-\Vert B\Vert I_n)A^\t\le 0$, which obviously holds for any
nonnegative definite matrix $B\in\R^{n\times n}$ and any $A\in\R^{m\times n}$.} 
\begin{eqnarray}\nonumber
&& \Phi^{-1}(u)(-\hat D^s(u))^{1\over 2}(-\hat D^s(u))^{-1\over 2}SU(u)U(u)^\t S^\t  
(-\hat D^s(u))^{-1\over 2}(-\hat D^s(u))^{1\over 2}(\Phi^{-1}(u))^\t \\
\lbl{integrand}
&\le& -\Vert F(u)\Vert\Phi(u)^{-1}\hat D^s(u)(\Phi(u)^{-1})^\t 
= {1\over 2}\Vert F(u)\Vert {d\over du}\{ \Phi(u)^{-1}(\Phi(u)^{-1})^\t\}
\end{eqnarray}
where 
$$
F(u):=(-\hat D^s(u))^{-1\over 2}SU(u)U^\t (u)S^\t (-\hat D^s(u))^{-1\over 2}.
$$
By taking into account (\ref{integrand}), from (\ref{covy}) we
have
\be\lbl{almost}
\tr\cov y_t\le{\sigma^2\over 2} \sup_{u\in[0,t]}\{\Vert F(u)\Vert\} 
\tr\{ I-\Phi(t)^\t\Phi(t)\}\le {\sigma^2n\over 2}
\sup_{u\in[0,t]}\Vert F(u)\Vert .
\ee
Further \footnotemark[1],
\begin{eqnarray}\nonumber
\Vert F(u)\Vert&=&
\Vert (-\hat D^s(u))^{-1\over 2}SU(u)U^\t (u)S^\t (-\hat D^s(u))^{-1\over 2}\Vert\\
&\le&
\lbl{there}
\Vert U(u)U^\t (u)\Vert \Vert SS^\t\Vert \Vert (-\hat D^s(u))^{-1}\Vert\\
\nonumber
& \le&  \alpha^{-1}\Vert U(u)U^\t (u)\Vert,
\end{eqnarray}
since $\Vert SS^\t \Vert = 1$ (cf. (\ref{orthogonality})).
Estimate (\ref{main}) follows from (\ref{Ey}), (\ref{almost}), and 
(\ref{there}).\\
\qed

Theorem~\ref{robust} describes a class of stochastically stable
time-dependent CPs. Because much of the previous work focused
on CPs with constant coefficients, we study them separately.
To this end, we consider
\be\lbl{const-perturbed}
\dot x=Dx +\sigma\dot w,
\ee
and the corresponding reduced system for $y=Sx$
\be\lbl{const-red}
\dot y =\hat Dy +\sigma S\dot w.
\ee
In (\ref{const-perturbed}), we set $U(t)=I_n$ to simplify notation.

\begin{thm}\lbl{first}
Suppose that $D\in\mathcal{D}$ and $\alpha$ has the same
meaning as in (\ref{alpha}). Then for any $0<\epsilon<\alpha$, 
there exists a positive constant $C_2$ such that
\be\lbl{first-estimate}
\lim_{t\to\infty} \E P_{\one^\perp}x(t)=0\quad\mbox{and}\quad
\E|P_{\one^\perp}x(t)|^2 
 \le {C_2\sigma^2 n \over \alpha-\epsilon}.
\ee
\end{thm}
\pf
Since $D\in\mathcal{D}$, for any $\epsilon\in (0,\alpha)$, there exists 
a change of coordinates in $\R^n$, $x=Q_\epsilon \tilde x,$ $Q_\epsilon\in\R^{n\times n}$,
such that $\tilde D_\epsilon = Q^{-1}_\epsilon DQ_\epsilon\in D_{\alpha-\epsilon}$
(cf. \cite{Arnold-ODEs}). By Theorem~\ref{robust}, solutions of 
$$
\dot{\tilde x}=\tl D_\epsilon\tilde x+\sigma Q^{-1}_\epsilon\dot w
$$
satisfy (\ref{main}). Thus, (\ref{first-estimate}) holds with 
for some $C_2>0$ possibly depending on $\epsilon$.\\
\qed

The estimate of $\E |P_{\one^\perp}x(t)|^2$ in (\ref{first-estimate})
characterizes the dispersion of the trajectories of the stochastically
forced CP (\ref{const-perturbed}) around the consensus subspace. 
$\E |P_{\one^\perp}x(t)|^2$ can be viewed as a measure of stochastic stability
of consensus subspace. In (\ref{first-estimate}), the upper bound for 
$\E |P_{\one^\perp}x(t)|^2$ is given in terms of the leading nonzero 
eigenvalue of $D$. If $D$ is normal then precise asymptotic values
for $\cov Sx(t)$ and  $\E |P_{\one^\perp} x(t)|^2$ are available.
The stability analysis of (\ref{const-perturbed}) with normal $D$ 
is important for understanding the properties of CPs on undirected graphs,
which we study in the next section.

\begin{thm}\lbl{variability}
Suppose $D\in\mathcal{D}$ is normal. Denote the EVs of $D$ 
by $\lambda_1,\lambda_2,\dots\lambda_n$, where $\lambda_1=0$
is a simple EV. Let $\hat D$ be a pseudo-similar matrix to $D$
via $S$ (cf. (\ref{rescaleS})).
Then for any deterministic initial condition $x(0)\in\R^n$,
the trajectory of (\ref{const-perturbed}) is a Gaussian random
process with the following asymptotic properties 
\begin{eqnarray}
\lbl{asympt-mean}
\lim_{t\to\infty} |\E P_{\one^\perp}x(t)|&=&0,\\
\lbl{asympt-cov}
\lim_{t\to\infty} \cov Sx(t)&=&2^{-1}\sigma^2 \hat D^{-1},\\
\lbl{asympt-norm}
\lim_{t\to\infty} \E |P_{\one^\perp} x(t)|^2&=&
2^{-1}\sigma^2 \sum_{i=2}^n (\mathrm{Re} \lambda_i)^{-1}.
\end{eqnarray}
\end{thm}
\pf 
By the observation in Remark~\ref{note},
it sufficient to prove the relations (\ref{asympt-mean}) and 
(\ref{asympt-norm}) with $P_{\one^\perp} x(t)$ replaced by $y(t)=Sx(t)$.

The solution of the reduced equation (\ref{const-red}) with a
deterministic initial condition  is a Gaussian process (cf. \cite{KS})
\be\lbl{variation}
y(t)=e^{t\hat D}y_0+\sigma\int_0^t e^{(t-u)\hat D}Sdw(u).
\ee
From (\ref{linear})
specialized to solutions of (\ref{const-red}), we have
\begin{eqnarray} \lbl{mean-y}
\E y(t) &=& e^{t\hat D} y_0\to 0, \mbox{as} \;t\to\infty,\\
\lbl{cov}
V(t)=\cov y(t) &=& \sigma^2 \int_0^t e^{(t-u)\hat D} SS^\t e^{(t-u)\hat D^\t }du.
\end{eqnarray} 
Since $SS^\t=I_{n-1}$ and $\hat D$ is a stable normal matrix, from (\ref{cov}) we have
\be\lbl{simply-integrate}
V(t)=\int_0^t e^{2u\hat D^s} du \to {\sigma^2\over 2} (\hat D^s)^{-1}, \; 
t\to \infty,
\ee
where $\hat D^s=2^{-1}(\hat D+\hat D^\t)$ stands for the symmetric part of $\hat D$.
By taking into account (\ref{mean-y}), we have
$$
\lim_{t\to\infty}\E|y(t)|^2= \lim_{t\to\infty}\tr V(t)= {\sigma^2\over 2} \tr \hat D^s=
 {\sigma^2\over 2} \sum_{i=2}^n (\mathrm{Re} \lambda_i)^{-1}.
$$
\qed
\begin{rem} Estimate (\ref{asympt-norm}) was derived in \cite{Leonard10}
for CPs with positive weights. A similar estimate was obtained in 
\cite{medvedev09} in the context of analysis of a neuronal network.
\end{rem}

\section{CPs on undirected graphs}\lbl{undirected}
\setcounter{equation}{0}
In this section, we apply the results of the previous section 
to CPs on undirected graphs. The analysis reveals the contribution 
of the network topology to stability of CPs. In particular, 
we show that the dimension and the structure of the cycle subspace associated 
with the graph of the network are important for stability. 
The former quantity, the first Betti number of the graph, is a 
topological invariant of the graph of the network.  
 
\subsection{Graph-theoretic preliminaries}
We start by reviewing certain basic algebraic constructions used
in the analysis of graphs (cf. \cite{Biggs}).
Throughout this subsection, we assume that $G=( V(G),E(G))$ 
is a connected graph with $n$ vertices and $m$ edges:
$$
V(G)=\left\{ v_1, v_2, \dots, v_n\right\}\quad\mbox{and}\quad
E(G)=\left\{ e_1, e_2, \dots, e_m\right\}.
$$
The vertex space, $C_0(G)$, and the edge space, $C_1(G)$, are the
finite-dimensional vector spaces of real-valued functions on 
$V(G)$ and $E(G)$, respectively.

We fix an orientation on $G$ by assigning positive and negative ends for
each edge in $E(G)$. 
The matrix of the coboundary operator $H:~C_1(G)\to C_0(G)$ with 
respect to the standard bases in $C_{0}(G)$ and $C_{1}(G)$ is defined by 
\be\lbl{incidence}
H=(h_{ij})\in \R^{m\times n},\quad
h_{ij}=\left\{ \begin{array}{cl}
1, & v_j\;\mbox{ is a positive end of}\; e_i,\\
-1, & v_j\;\mbox{ is a negative end of}\; e_i,\\
0, &\;\mbox{otherwise}.
\end{array}
\right.
\ee
The Laplacian of $G$ is expressed in terms of the coboundary matrix
\be\lbl{Lap}
L=H^\t H.
\ee

By $\tilde G=( V(\tl G),E(\tl G))\subset G$ we denote a spanning tree of $G$,
a connected subgraph of $G$ such that
$$
\left|V(\tl G)\right|=n \quad\mbox{and}\quad \left|E(\tl G)\right|=n-1.
$$
$\tilde G$ contains no cycles. Without loss of generality, we assume
that 
\be\lbl{span-tree}
E(\tilde G)=\{e_1, e_2, \dots, e_{n-1}\}.
\ee

A cycle $O$ of length $|O|=k$ is a cyclic subgraph of $G$:
$$
O=( V(O), E(O))\subset G:\quad
V(O)=\{v_{i_1},v_{i_2},\dots,v_{i_k}\},\quad 
E(O)=\{(v_{i_1},v_{i_2}), (v_{i_2},v_{i_3})\dots,(v_{i_k},v_{i_1})\},\;
$$
for some $k$ distinct integers $(i_1, i_2,\dots, i_k)\in [n]^k$.
Two cyclic ordering of the vertices of $O$ induce two
possible orientations. Suppose the orientation of $O$
has been fixed. We will refer to it as the cycle orientation.
For each $e_i\in E(O)$, we thus have two orientations: one induced
by the orientation of $G$ and the other induced by the cycle orientation.
We define $\xi(O)=(\xi_1,\xi_2,\dots, \xi_m)^\t \in C_1(G)$
such that
\be\lbl{define-xi}
\xi_i=\left\{\begin{array}{cl}
1, & \mbox{if $e_i\in E(O)$ and the two orientations of $e_i$ coincide},\\ 
-1, & \mbox{if $e_i\in E(O)$ and the two orientations of $e_i$ differ},\\ 
0, &\mbox{otherwise}.
\end{array}
\right.
\ee
To each cycle $O$ of $G$, there corresponds $\xi(O)\in C_1(G)$.
All such vectors span the cycle subspace of $G$,
$\mathrm{Cyc}(G)\subset C_1(G)$. The cycle subspace coincides
with the kernel of the incidence mapping $H^\t $
\be\lbl{kerH}
\mathrm{Cyc}(G)=\ker~H^\t ,
\ee
and its dimension is equal to the corank of $G$
\be\lbl{define-c}
c=m-n+1.
\ee

To each edge $e_{n-1+k}, k\in [c],$ not belonging to the 
spanning tree $\tl G$, there corresponds a unique cycle $O_k$
such that 
$$
e_{n-1+k}\in E(O_k)\subset E(\tl G)\cup\{e_k\}.
$$
We orient cycles $O_k, k\in [c],$ in such a way that the orientations of 
$e_{n-1+k}$ as an edge of $G$ and that of $O_k$ coincide.
The vectors 
\be
\xi^k=\xi(O_k),\; k\in [c],
\ee
form a basis in $\mathrm{Cyc}~(G)$. We will refer to cycles $O_k, k\in [c]$
as the fundamental cycles of $G$.

Define $c\times m$ matrix
$$
Z=\begin{pmatrix} {\xi^1}^\t  \\{\xi^2}^\t \\\dots \\{\xi^c}^\t  \end{pmatrix}.
$$
By construction, $Z$ has the following block structure,
\be\lbl{blocks-of-Z}
Z=(Q\; I_c),
\ee
where 
$I_c$ is the $c\times c$ identity matrix. Matrix $Q=(q_{kl})\in\R^{c\times (n-1)}$
contains the coefficients of expansions of $\vec{e}_{n-1+k}, k\in[c]$
with respect to $\{ \vec{e_i},\; i\in [n-1]\}$
\be\lbl{expand-cycle}
\vec{e}_{n-1+k}=-\sum_{l=1}^{n-1} q_{kl}\vec{e_l},\; q_{kl}\in\{0, \pm 1\}.
\ee
Here, $\vec{e_i},\; i\in[m]$ denote the oriented edges of $G$. This motivates the 
following definition.
\begin{df}
\begin{enumerate}
\item Matrix $Q\in\R^{c\times (n-1)}$ is called a cycle incidence matrix of
$G$.
\item Matrix $L_c(G)=QQ^\t \in \R^{c\times c}$ is  called a cycle Laplacian
of $G$.
\end{enumerate}
\end{df}   

The following properties of the cycle Laplacian follow
from its definition. 
\begin{description}
\item[A] The spectrum of $L_c(G)$ does not depend on the choice of orientation used 
in the  construction of $Q$. Indeed, if $Q_1$ and $Q_2$ are two cycle incidence matrices 
of $G$ corresponding to two distinct orientations then $Q_2=P_c Q_1 P_t$, where diagonal matrices 
$$
P_t=\mbox{diag}\{p_1, p_2,\dots, p_{n-1}\}\in \R^{(n-1)\times (n-1)}
\quad\mbox{and}\quad
P_c=\mbox{diag}\{p_n, p_{n+1},\dots, p_m\}\in \R^{c\times c},
$$ 
If the two orientations of $G$ yield the same orientation for edge $e_i\in E(G), i\in [m]$,
then $p_i=1$; and $p_i=-1$, otherwise. 
Thus,
$$
Q_2 Q_2^\t=P_cQ_1 P_tP_t^\t Q_1^\t P_c^\t = P_c Q_1 Q_1^\t P_c^\t.
$$
The spectra of $Q_2 Q_2^\t$ and $Q_1 Q_1^\t$ coincide, because $P_c$ is an orthogonal matrix.
\item[B] 
\be\lbl{dot}
(L_c(G))_{ij}=
\langle\mathrm{Row}_i(L_c(G)),\mathrm{Row}_j(L_c(G))\rangle=
\left\{ \begin{array}{cc}
|O_i|-1, & i=j,\\
\pm|O_i\cap O_j|,& i\ne j.
\end{array}
\right.
\ee
\item[C]
If the cycles are disjoint, then 
assuming that that $O_k, k\in[c],$
are ordered by their sizes, we have
\be\lbl{disjoint-EVs}
\lambda_k(I_c+QQ^\t )=|O_k|,\; k\in [c].
\ee
\end{description}

The cycle incidence matrix provides a convenient partition
of the coboundary matrix.
\begin{lem}\lbl{partition}
Let $G=(V(G),E(G))$ be a connected graph of positive 
corank $c$.
Then using the notation of this section,
the coboundary matrix 
\be\lbl{part}
H=\begin{pmatrix} I_{n-1}\\ -Q\end{pmatrix} \tl H,
\ee
where $\tl H$ is the coboundary matrix of a spanning tree
$\tl G$, and $Q$ is the corresponding cycle incidence matrix.
\end{lem}
\pf 
Since $\mathrm{rank} H=\mathrm{rank} \tl H=n-1$,
there is a unique $B\in\R^{c\times (n-1)}$ such that 
\be\lbl{lin-ind}
H=\begin{pmatrix}
\tl H \\ B\tl H
\end{pmatrix}.
\ee
Using (\ref{lin-ind}) and (\ref{blocks-of-Z}),
we obtain
$$
ZH=0\;\Rightarrow\; Q\tl H+B\tl H=0\;\Rightarrow\; B=-Q.
$$  
\qed

\subsection{Stability analysis}
We are now in a position to apply the results of Section~\ref{robustness}
to study CPs on undirected graphs.
Let  $G=(V(G), E(G))$ be a connected undirected graph.
Since the interactions between the agents are symmetric, the 
coupling matrix in (\ref{const-perturbed}) can be 
rewritten in terms the coboundary matrix of $G$ and the 
conductance matrix $C=\mathrm{diag} (c_1,c_2,\dots,c_m)$
\be\lbl{weight-Lap}
D=-H^\t C H.
\ee
The conductances $c_i\in\R, i\in [m]$ are assigned to all edges 
of $G$. If the network is simple, $C=I_m$.

Let $\tilde G$ be a spanning tree of $G$. 
We continue to assume that the edges of $G$ are ordered
so that (\ref{span-tree}) holds. Let $\tilde H$, $Q$,
and $C=\diag (C_1,C_2),$ $C_1\in\R^{(n-1)\times (n-1)}$, 
$C_2\in\R^{c\times c}$,  be the coboundary,  cycle
incidence matrix and conductance matrices 
corresponding to the spanning tree $\tl G$, respectively.
Using (\ref{part}),  we recast the coupling matrix  as
\be\lbl{recast-coup}
D=-{\tl H}^\t (C_1+Q^\t C_2Q)\tl H.
\ee

To form the reduced equation, we let
$$
S=(\tl H {\tl H}^\t)^{-{1\over 2}}\tl H \quad \mbox{and}\quad y=Sx.
$$
Then
\be\lbl{new-red}
\dot y=\hat D y+ \sigma S\dot w,
\ee
where
\be\lbl{rewrite}
\hat D=-(\tl H {\tl H}^\t)^{1\over 2} (C_1+Q^\t C_2Q)(\tl H {\tl H}^\t)^{1\over 2}.
\ee
Both $D$ and $\hat D$ are symmetric matrices. By Lemma~\ref{spectra},
the eigenspaces of $\hat D$ and $D$ are related via $S$. The EVs
of $\hat D$ are the same as those of $D$ except for a simple
zero EV $\lambda_1=0$ corresponding to the constant eigenvector 
$e$.

A CP with positive weights $c_i$'s is  always convergent, as can be
easily seen from (\ref{rewrite}). For a general case, we have the 
following necessary and sufficient condition for convergence
of the CP on an undirected graph.

\begin{thm}\lbl{undirected-conv}
The CP (\ref{constant}) with matrix (\ref{weight-Lap}) is convergent if 
and only if matrix $C_1+Q^\t C_2Q$ is positive definite
for some spanning tree $\tl G$. 
\end{thm}
\pf By (\ref{rewrite}), $\hat D$ is a stable matrix if and only if
$C_1+Q^\t C_2Q$ is positive definite.\\
\qed

If $D\in\mathcal{D}$ stochastic stability of the CP (\ref{const-perturbed}$)_0$  and 
(\ref{weight-Lap}) is guaranteed by Theorem~\ref{variability}. In particular, (\ref{asympt-norm})
characterizes the dispersion of trajectories around the consensus subspace.
Theorems~\ref{exhaustive} and \ref{variability} provide explicit formulas
for the rate of convergence and the degree of stability of convergent CPs 
on undirected graphs. Specifically, let $\mathcal{N}=(G, C)$ be a network corresponding
to (\ref{const-perturbed}) with $D\in\mathcal{D}$ (cf. (\ref{recast-coup})).
Let 
$$
0=\lambda_1<\lambda_2\le\dots\le\lambda_n
$$
denote the EVs of $-D$ and define
\be\lbl{alpha-and-rho}
\alpha(\mathcal{N})= \lambda_2 \quad\mbox{and}\quad 
\rho(\mathcal{N})=\sum_{i=2}^n {1\over\lambda_i}.
\ee
Formulas in (\ref{alpha-and-rho}) generalize algebraic connectivity and
(up to a scaling factor) total effective resistance of a simple graph to 
weighted networks corresponding
to convergent CPs on undirected graphs.
By replacing the EVs of $D$ by those of $D^s$ in (\ref{alpha-and-rho}),
the definitions of $\alpha(\mathcal{N})$ and $\rho(\mathcal{N})$ can be 
extended to convergent CPs with normal coupling  matrices.
For simple networks, there are many results relating algebraic connectivity and
total effective resistance and the structure of the graph 
(cf. \cite{Fiedler73, Allon86, Boyd08,Chung97, Cioaba10,Gutman03, Hoory06, Nil91}).
Theorems~\ref{exhaustive} and ~\ref{variability} link 
structural properties of the network and dynamical performance of the CPs.

In conclusion of this section, we explore some implications of (\ref{rewrite})
for stability of (\ref{const-perturbed}). To make the role of the network topology in shaping
stability properties of the system more transparent, in the remainder of this 
paper, we consider simple networks, i.e., $C=I_m$. In the context of stability of
CPs, the lower bounds for $\alpha(\mathcal{N})$ 
and the upper bounds for $\rho(\mathcal{N})$ are important.  
\begin{lem}\lbl{one-more}
Let $\tl H$ stand for the coboundary matrix of a spanning tree $\tl G$ 
of undirected graph $G$ and let $Q$ be the corresponding  cycle incidence 
matrix. Then
\begin{description}
\item[A]
$$
\alpha(G)\ge 
\alpha (\tl G) \lambda_1\left(I_{n-1}+Q^\t Q\right),
$$
\item[B]
$$
\rho(\mathcal{N})\le 
\min\left\{{n-1\over\lambda_1\left(I_{n-1}+Q^\t Q\right)},  
{\tr (I_{n-1}+Q^\t Q)^{-1}\over\alpha(\tl G)} \right\}
\le (n-1)\min\left\{1, {1\over\alpha(\tl G)}\right\},
$$
where $\alpha(G)$ 
stands for the algebraic connectivity of $G$, and 
$\lambda_1\left(I_{n-1}+Q^\t Q\right)$ denotes the 
smallest EV of the positive definite matrix  $I_{n-1}+Q^\t Q$.
\end{description}
\end{lem} 
\pf Since $\mathcal{N}$ is a simple network, the coupling matrix
$D$ taken with a negative sign is the Laplacian of $G$, $L=H^\t H$.
Likewise, $\tl L= {\tl H}^\t \tl H$ is a Laplacian of $\tl G$.
Let
$$
0=\lambda_1(G)<\lambda_2(G)\le\dots\le\lambda_n(G)
$$
denote the EVs of $L$. Below we use the same notation to denote  the EVs
of other positive definite matrices, e.g.,  $\tl L$ and $I_{n-1}+Q^\t Q$. 

By Lemma~\ref{spectra}, the second EV of $G$, $\alpha(G)$, coincides
with the smallest  EV of 
$$
\hat L=-\hat D
=({\tl H}{\tl H}^\t)^{1\over 2}(I_{n-1}+Q^\t Q)({\tl H}{\tl H}^\t)^{1\over 2}.
$$
Below, we will use the following observations.
\begin{description}
\item[a]
The sets of nonzero EVs of two symmetric matrices $\tl H {\tl H}^\t$
and $\tilde L=\tl H^\t {\tl H}$ coincide. Since the former is a full rank matrix,
the spectrum of $\tl H {\tl H}^\t$ consists of nonzero EVs
of $\tl G$. In particular,
\be\lbl{equalEVs}
\lambda_1(\tl H{\tl H}^\t)=\alpha(\tl G).
\ee
\item[b] The EVs of $\hat L$ and those of $(\tl H{\tl H}^\t) (I_{n-1}+Q^\t Q)$
coincide.
\end{description}
Using the variational characterization of the EVs of symmetric matrices (cf. \cite{HJ99})
and the observations \textbf{a} and \textbf{b} above, we have
$$
\lambda_1(\hat L)=\lambda_1\left( (\tilde H {\tilde H}^\t)(I_{n-1}+Q^\t Q)\right)\ge
\lambda_1(\tl H{\tl H}^\t)\lambda_1(I_{n-1}+Q^\t Q)=\alpha(\tl G)\lambda_1(I_{n-1}+Q^\t Q). 
$$
Hence,
$$
\alpha(G)=\lambda_1(\hat L)\ge \alpha(\tl G)\lambda_1(I_{n-1}+Q^\t Q).
$$
Likewise\footnote[1]{Estimate (\ref{elementary}) uses the following inequality:
$\lambda_1 (A)\tr B\le \tr (AB)\le\lambda_n(A)\tr B$, which holds  for any symmetric 
$A\in\R^{n\times n}$ and nonnegative definite $B\in\R^{n\times n}$ (cf.~\cite{KA}). 
If $A=\Lambda$ is diagonal, the double inequality above is obvious. The case of symmetric
$A=U\Lambda U^\t\;$ is reduced to the previous one via the following identity 
$\tr\{U \Lambda U^\t B\}=\tr\{\Lambda U^\t BU\}$.},
\begin{eqnarray}\nonumber
\rho(G)
&=&\tr (-{\hat D}^{-1})=\tr \left\{ (\tl H {\tl H}^\t)^{-1} (I_{n-1}+Q^\t Q)^{-1}\right\} 
\\ 
\lbl{elementary}
&\le& {\tr \{(\tl H{\tl H}^\t)^{-1} \} \over \lambda_1(I_{n-1}+Q^\t Q)}=
{ \rho(\tl G) \over \lambda_1(I_{n-1}+Q^\t Q)}={ n-1\over \lambda_1(I_{n-1}+Q^\t Q)}.
\end{eqnarray}
A symmetric argument yields
$$
\rho(G)\le {\tr (I_{n-1}+Q^\t Q)^{-1}\over \lambda_1(\tl H{\tl H}^\t)}=
{\tr (I_{n-1}+Q^\t Q)^{-1}\over\alpha(\tl G)}.
$$
\qed   

Lemma~\ref{one-more} shows respective contributions 
of the spectral properties of the spanning tree $\tl G$ and those of the cycle
subspace to the algebraic connectivity and effective resistance 
of $G$. In this respect, it is of interest to study the 
spectral properties of $I_{n-1}+Q^\t Q$, in particular, its smallest
EV and the trace. Another motivation for studying $I_{n-1}+Q^\t Q$
comes from the following lemma.

\begin{lem}\lbl{rephrase} 
Under assumptions of Lemma~\ref{one-more}, solutions of CP 
(\ref{const-perturbed}) satisfy
\be
\lbl{kappa}
\lim_{t\to\infty} \E |\tl H x(t)|^2 =
{\sigma^2\over 2}\kappa(G,\tl G),\quad \kappa(G,\tl G)=\tr (I_{n-1}+Q^\t Q)^{-1}.
\ee
\end{lem}
\pf 
The reduced equation for $y=\tl Hx$ has the following form
\be\lbl{redd}
\dot y=\hat Dy+\sigma\tl H \dot w,
\ee
where
\be\lbl{dhatt}
\hat D= \tl H D {\tl H}^+.
\ee
Using $D=-H^\t H$ and (\ref{part}), we rewrite (\ref{dhatt}) as follows
\be\lbl{dhatt1}
\hat D= \tl H D {\tl H}^+
= -\tl H H^\t H {\tl H}^\t (\tl H {\tl H}^\t)^{-1}=-\tl H {\tl H}^\t (I_{n-1}+Q^\t Q).
\ee
By applying the argument used in the proof of Theorem~\ref{variability}
to the reduced equation (\ref{redd}), we obtain
\be\lbl{new-lim}
\lim_{t\to\infty} \E |\tl H x(t)|^2 =
{\sigma^2 \over 2} \tr \{\tl H {\tl H}^\t (\hat D)^{-1} \}={\sigma^2\over 2} \kappa (G,\tl G).
\ee
The combination of (\ref{new-lim}) and (\ref{dhatt}) yields (\ref{kappa}).\\
\qed

The following lemma provides the graph-theoretic 
interpretation of $\kappa(G,\tl G)$.

\begin{lem}\lbl{cycles}
Let $G=(V(G), E(G)), \; |V(G)|=n,$ be a connected graph and 
$\tl G\subset G$ be a spanning tree of $G$.
\begin{description}
\item[A] If $G$ is a tree then
\be\lbl{n-1}
\kappa(G,\tl G)=n-1.
\ee
\item[B] Otherwise, denote the corank of $G$ by $c>0$ and let $\{O_k\}_{k=1}^c$ be the 
system of fundamental cycles corresponding to $\tl G$.
\begin{description}
\item[B.1] Denote 
\be\lbl{weight}
\mu={1\over n-1}\sum_{k=1}^c (|O_k|-1).
\ee
Then
\be\lbl{cool}
{1\over 1+\mu} \le {\kappa(G,\tl G)\over n-1}\le 1,
\ee
\item[B.2] If $0<c<n-1$ then 
\be\lbl{c-small}
1-{c\over n-1}\left(1-{1\over \delta}\right)\le {\kappa(G, \tilde G)\over n-1}
\le 1,
\ee
where
\be\lbl{interposition}
\delta=\max_{k\in [c]} \{|O_k| +\sum_{l\neq k} |O_k\cap O_l|\}.
\ee
\item[B.3]
If $O_k, k\in [c]$ are disjoint.
Then 
\be\lbl{disjoint}
{\kappa(G,\tl G)\over n-1}=
1-{c\over n-1}\left(1-{1\over c}\sum_{k=1}^c |O_k|^{-1} \right).
\ee
In particular,
$$
{\kappa(G,\tl G)\over n-1} \le
1-{c\over n-1}\left(1-{1\over \min_{k\in [c]} |O_k|} \right)
$$
and 
$$
{\kappa(G,\tl G)\over n-1}\ge 
1-{c\over n-1}\left(1-{c\over \sum_{k\in [c]} |O_k|} \right)
\ge 
1-{c\over n-1}\left(1-{1\over \max_{k\in [c]} |O_k|} \right).
$$
\end{description}
\end{description}
\end{lem}
\pf
\begin{description}
\item[A]
If $G$ is a tree then $Q=0$ and $\kappa(G,\tl G)=\tr I_{n-1}=n-1$.
\item[B.1]
Suppose $c>0$.
Let $\lambda_i,\;i\in [n]$, denote the EVs of $I_{n-1}+Q^\t Q$:
$$
1\le\lambda_1\le\lambda_2\le\dots\le\lambda_{n-1}
$$
By the arithmetic-harmonic means
inequality, we have
\be\lbl{harmonic}
\left({1\over n-1}\sum^{n-1}_{i=1}\lambda_i\right)^{-1}\le 
{1\over n-1}\sum^{n-1}_{i=1}{1\over \lambda_i}\le {1\over \min_{i\in [n-1]}\lambda_i}\le 1.
\ee
The double inequality in (\ref{cool}) follows from (\ref{harmonic})
by noting that 
$$
\kappa (G,\tl G)=\sum^{n-1}_{i=1}{1\over \lambda_i}
$$
and
$$
\sum^{n-1}_{i=1}\lambda_i=\tr(I_{n-1}+Q^\t  Q)=n-1+\tr QQ^\t =n-1+\sum^c_{k=1}\left(|O_k|-1\right).
$$
\item[B.2]
Since $\mathrm{rank}~Q^\t  Q=c<n-1$, by the interlacing theorem 
(cf. Theorem~4.3.4 \cite{HJ99}), we have
\be\lbl{n-c}
1\le\lambda_k(I_{n-1}+Q^\t  Q)\le \lambda_{k+c}(I_{n-1})=1,\; k\in [n-1-c].
\ee
For $k>n-1-c$, we use the Weyl's Theorem to obtain
\be\lbl{weyl}
1\le \lambda_k(I_{n-1}+Q^\t  Q)\le 1+\lambda_{n-1}(Q^\t  Q)=1+\lambda_{c}(QQ^\t ).
\ee
Using (\ref{dot}), by the Gershgorin's Theorem, we further have
\be\lbl{gersh}
1+\lambda_{c}(QQ^\t )\le \max_{k\in [c]} \{|O_k| +\sum_{l\neq k} |O_k\cap O_l|\}.
\ee
The combination of (\ref{weyl}) and (\ref{gersh}) yields
$$
\kappa(G,\tl G)=\sum_{k=1}^{n-1} {1\over\lambda_i}\ge n-1 +{c\over\delta}.
$$
\item[B.3]
Since each cycle $O_k$ contains at least two edges from the spanning
tree $\tl G$, then the number of disjoint cycles $c$ can not
exceed the integer part of $0.5(n-1)$. In particular, $c<n-1$.

By (\ref{dot}),
$$
QQ^\t =\mathrm{diag}(|O_1|-1,|O_2|-1,\dots,|O_c|-1),
$$  
because the cycles are disjoint. Further, the nonzero eigenvalues
of $Q^\t  Q$ and $QQ^\t $ coincide. Thus,
\be\lbl{sp}
\lambda_k\left((I_{n-1}+Q^\t  Q)^{-1}\right)=\left\{
\begin{array}{cc}  
1, & k\in [n-1-c], \\
| O_{k+c+1-n}|^{-1}, & n-c\le k\le n-1,
\end{array}
\right.
\ee
By plugging (\ref{sp}) in (\ref{kappa}), we obtain (\ref{disjoint}).
\end{description}
\qed
\begin{rem}
Estimates of $\kappa(G,\tl G)$ in Lemma~\ref{cycles}, combined with 
the estimates in Lemmas~\ref{one-more} and \ref{rephrase}, show
how stochastic stability of CPs depends on the geometric properties 
of the cycle subspace associate with the graph, such as 
the first Betti number (cf. (\ref{c-small})) and
the length and the mutual position the fundamental cycles 
(cf. (\ref{weight}), (\ref{cool}), (\ref{interposition}), (\ref{disjoint})).
In particular, from (\ref{kappa}) and the estimates in Lemma~\ref{cycles}
one can see how the changes of the graph of the network, which do not affect a 
spanning tree, impact stochastic stability of the CP. Likewise, by combining
the statements in Lemma~\ref{cycles} with the following estimate
of the total effective resistance (cf. Lemma~\ref{one-more})
\be\lbl{very}
\rho(\mathcal{N})\le \alpha(\tl G)\kappa(G,\tl G),
\ee 
one can see how the properties of the spanning tree and the corresponding
fundamental cycles contribute to the stochastic stability of the CP.
\end{rem}
\begin{figure}
\begin{center}
{\bf a}\;\epsfig{figure=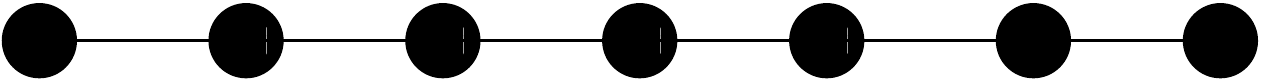, height=0.15in, width=1.8in, angle=0}
{\bf b}\;\epsfig{figure=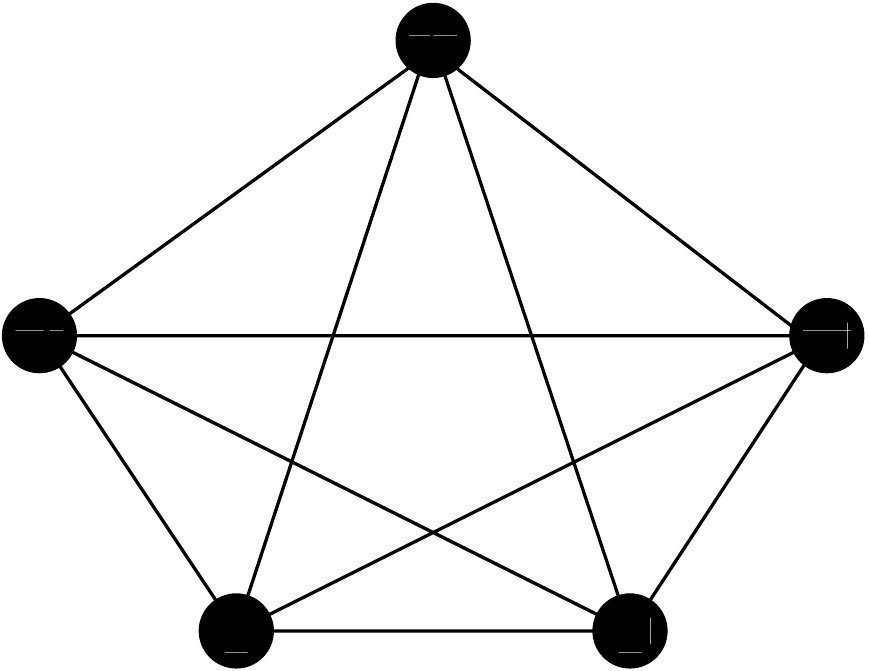, height=0.75in, width=1.8in, angle=0}
{\bf c}\;\epsfig{figure=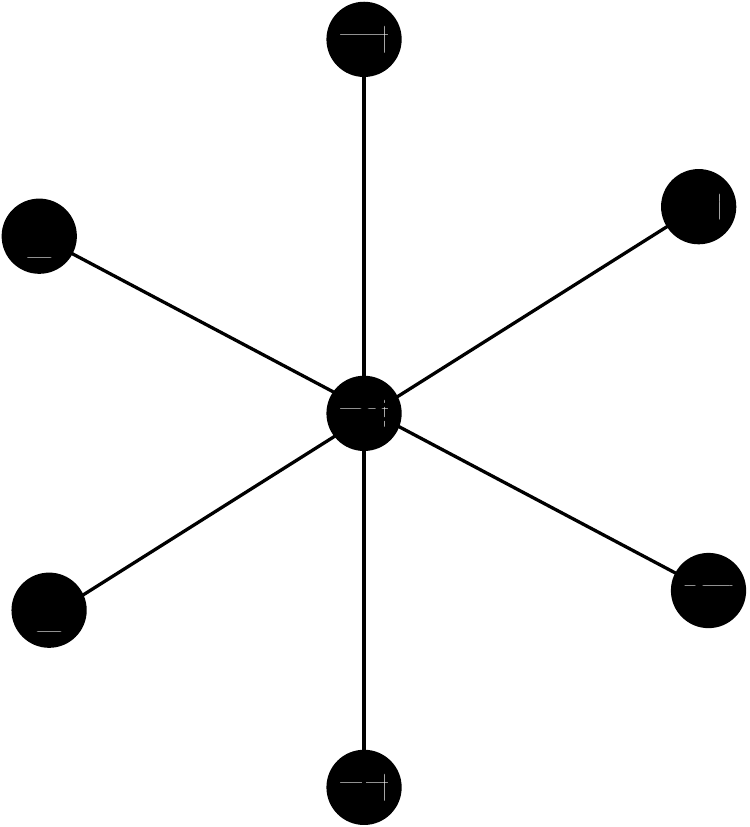, height=0.85in, width=1.8in, angle=0}
\end{center}
\caption{Examples of graphs used in the text:
(a) a path, (b) a complete graph (see Example~\ref{pair} for discussion
of the properties of graphs in (a) and (b)), and
(c) a star.
}\label{f.0}
\end{figure}

\section{Network connectivity and performance of CPs}\lbl{connectivity} 
\setcounter{equation}{0}

In the previous section, we derived several quantitative estimates 
characterizing convergence and stochastic stability of CPs.
In this section, we discuss two examples illustrating how different structural 
features
of the graph shape the dynamical properties of CPs. In the first pair of examples,
we consider graphs of extreme degrees: $2$ vs. $n$. In the second example, we take
two networks of equal degree but with disparate connectivity patterns: random vs. 
symmetric. These examples show that both the degree of the network and its 
connectivity are important.
   

\begin{ex}\lbl{pair}
Consider two simple networks supported by  a path, $P_n$, and  by a complete graph, $K_n$
(Fig.~\ref{f.0} a and b). The coupling matrices of the corresponding CPs are given by
\be\lbl{path-complete}
D_p=
\left(\begin{array}{cccccc}
-1 & 1 & 0 & \dots &0 &0\\
1& -2& 1& \dots& 0& 0\\
\dots& \dots& \dots& \dots &\dots& \dots \\
0& 0& 0& \dots& 1& -1
\end{array}
\right) \; \mbox{and}\;
D_c=
\left(\begin{array}{cccc}
-n+1 & 1 &  \dots  & 1\\
1& -n+1&  \dots & 1\\
\dots& \dots&\dots& \dots \\
1& 1&  \dots& -n+1
\end{array}
\right).
\ee
\end{ex}
The nonzero EVs of $P_n$ and $K_n$ are given by 
\be\lbl{EVs:path-complete}
\lambda_{i+1}(P_n)=4\sin^2\left({\pi i\over 2n}\right)\quad\mbox{and}\quad
\lambda_{i+1}(K_n)=n,\; i=1,2,\dots, n-1.
\ee
Thus,
\be\lbl{compare}
\alpha (P_n)=4\sin^2\left({\pi\over 2n}\right),\; 
\alpha (K_n)=n,\;\mbox{and}\; \rho(K_n)=1-n^{-1}.
\ee
To compute compute $\rho(P_n)$, we use the formula for the total effective resistance
of a tree (cf. (5), \cite{Boyd08})
\be\lbl{rho-path}
\rho(P_n)=n^{-1}\sum_{i=1}^{n-1} i(n-i)={1\over 6} (n^2-1).
\ee
Equation (\ref{compare}) shows that for $n\gg 1$, the convergence rate of the CP based on
the complete graph is much larger than that based on the path:
\be\lbl{complete-vs-path}
\alpha(K_n)=n \gg O(n^{-2})=\alpha(P_n).
\ee
One may be inclined to attribute the disparity in the convergence
rates to the fact that the degrees of the underlying graphs (and,
therefore, the total number of edges) differ substantially.
To see to what extent the difference of the total number of the
edges or, in electrical terms, the amount of  wire 
used in the corresponding electrical circuits, can account for the mismatch
in the rates of convergence, we scale the coupling matrices in 
Example~\ref{pair} by the degrees of the corresponding graphs:
$$
\tilde D_p={1\over 2}D_p \quad\mbox{and}\quad 
\tilde D_c={1\over n-1}D_c.
$$
The algebraic connectivities of the rescaled networks
are still far apart:
\be\lbl{scale_path}
\lambda_2(\tilde D_c)=1+O(n^{-1})\gg O(n^{-2})=\lambda_2(\tilde D_p).
\ee
This shows that the different values of $\alpha$ in (\ref{complete-vs-path})
reflect the distinct patterns of connectivity of these networks.

\begin{rem}
Explicit formulas for the EVs of the graph Laplacian 
similar to those that we used for the complete graph and the 
path are available for a few other canonical coupling 
architectures such as a cycle, a star, an $m-$dimensional lattice 
(see, e.g.,  \S 4.4 \cite{Fiedler73}). Explicit examples of graphs
with known EVs can be used for developing intuition for 
how the structural properties of the graphs translate to
the dynamical properties of the corresponding CPs.
\end{rem}

Equation (\ref{complete-vs-path}) shows that the rate of convergence of CPs
based on local nearest-neighbor interactions decreases rapidly 
when the network size grows. Therefore, this network architecture
is very inefficient for coordination of large groups of agents. 
The following estimate   
shows that very slow convergence of the CPs based on a path 
for $n\gg 1$ is not specific to this particular network topology,
but is typical for networks with regular connectivity.
For graph $G_n$ of degree $d$ on $n$ vertices, the following
inequality holds \cite{AM85}
$$
\alpha(G_n)\le 2d \left[{2\log_2n\over \mbox{diam} (G_n)}\right]^2.
$$
This means that if the diameter of $G_n$ grows faster than $\log_2 n$ 
(as in the case of a path or any lattice), the algebraic connectivity 
and, therefore, the convergence rate of the CP goes to zero as $n\to\infty$.  
Therefore, regular network topologies such as lattices result in poor 
performance of CPs. In contrast, below we show that a random network 
with high probability has a much better (in fact, nearly optimal) 
rate of convergence.

The algebraic connectivity of the (rescaled) complete graph does not
decrease as the size of the graph goes to infinity (cf. (\ref{scale_path})). 
There is another canonical network architecture, a star (see Fig.~\ref{f.0}c),
whose algebraic connectivity
remains unaffected by increasing the size of the network: 
$$
\lambda_2 (D_s)=1.
$$ 
However, both the complete graph and the star have disadvantages
from the CP design viewpoint. CPs based on the complete graph are
expensive, because they require $O(n^2)$ interconnections. The star
uses only $n-1$ edges, but the performance of the entire network 
critically depends on a single node, a hub, that connects to all other nodes. 
In addition, update of the information state of the hub  
requires simultaneous knowledge of the states of all other agents in the
network.  Therefore, neither the complete graph nor the star 
can be used for distributed consensus algorithms.
 
\begin{figure}
\begin{center}
{\bf a}\quad\epsfig{figure=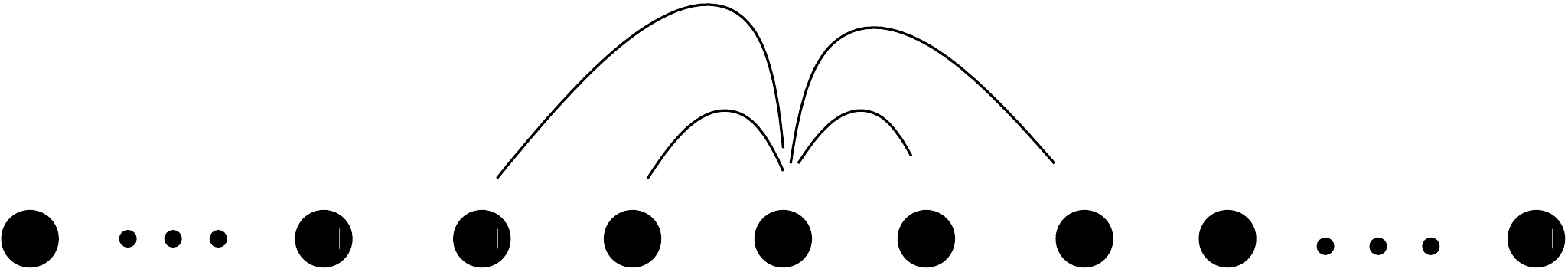, height=0.35in, width=2.5in, angle=0}
\qquad
{\bf b}\quad\epsfig{figure=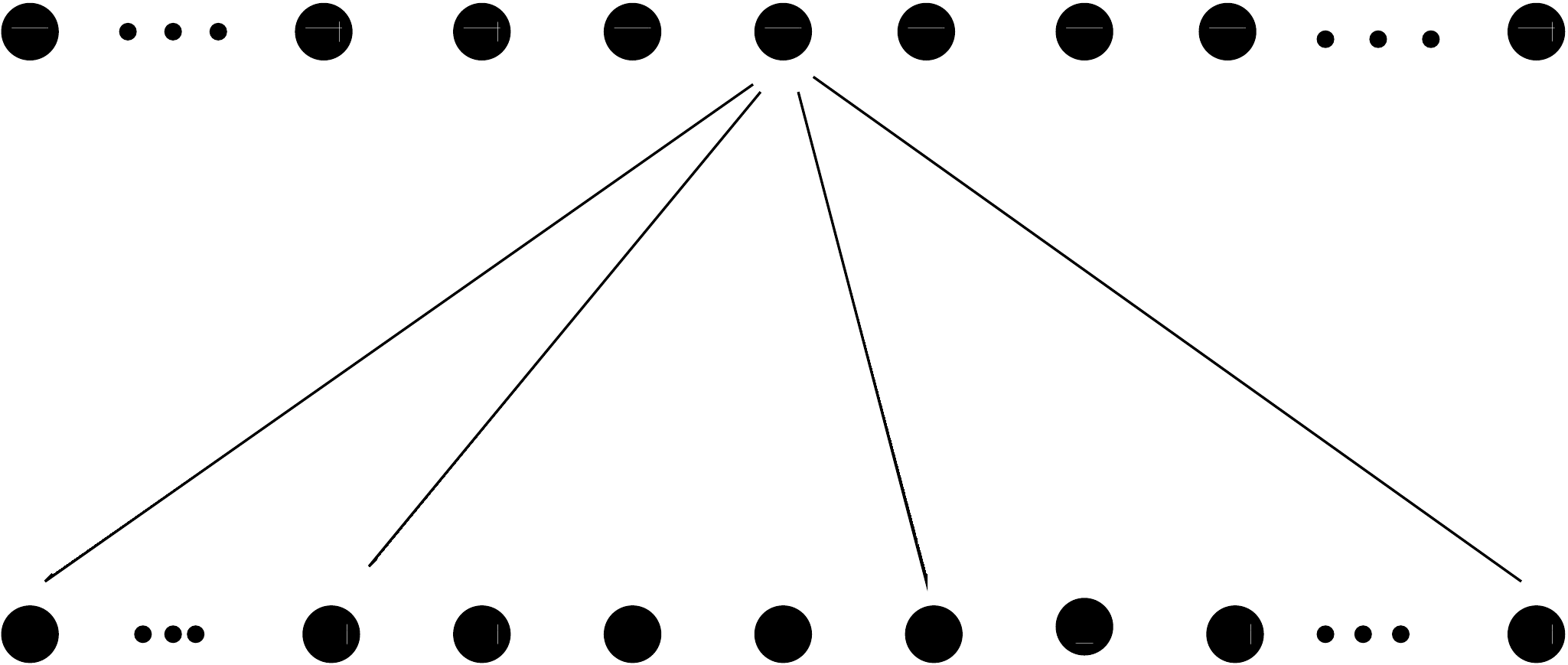, height=1.0in, width=2.5in, angle=0}
\end{center}
\caption{Schematic representation of the local connections
in  graphs used in Example~\ref{path_and_bipart}:
(a) a cycle of degree $4$ and (b) a bipartite graph of degree $4$ with
random connections.
}\label{f.1}
\end{figure}

Ideally, one would like to have a family of sparse graphs 
that behaves like that of complete graphs in the sense that the algebraic 
connectivity remains bounded from zero uniformly:
$$
\alpha (G_n)\ge \bar\alpha>0,\quad n\in\N.
$$
Moreover, the greater the value of $\bar\alpha$ the better the convergence
of the corresponding CPs.  Such graphs are called (spectral) expanders 
\cite{Hoory06,Sar04}. Expanders can be used for producing
CPs with a guaranteed rate of convergence regardless of the
size of the network. There are several known explicit constructions of 
expanders including celebrated Ramanujan graphs \cite{Margulis88, LPS88} 
(see \cite{Hoory06} for an excellent review of the theory and applications 
of expanders). In addition, random graphs are very good expanders. 
To explain this important property of random graphs, let us consider a family of graphs 
$\{G_n\}$ on $n$ vertices of fixed degree $d\ge 3$, i.e., $G_n$ is an 
$(n,d)-$graph. The following theorem due to Alon and Boppana yields 
an (asymptotic in $n$) upper bound on $\alpha(G_n)$.
\begin{thm}\cite{Nil91}
For any $(n,d)-$graph $G_n$, $d\ge 3$ and $n\gg 1$,
\be\lbl{AB}
\alpha(G_n)\le g(d)+o_n(1), \; g(d):=d-2\sqrt{d-1}>0.
\ee
\end{thm}
Therefore, for large $n$, $\alpha(G_n)$ can not exceed $g(d)$ 
more than by a small margin. The following theorem of Friedman
shows that for a random $(n,d)$-graph $G_n$, 
$\alpha(G_n)$ is nearly optimal with high probability.
\begin{thm}\lbl{random_graph}\cite{Fri08}
For every $\epsilon>0$,
\be\lbl{Friedman}
\mathsf{Prob}\left\{ \alpha(G_n)\ge g(d)-\epsilon\right\}=1-o_n(1),
\ee
where $\{G_n\}$ is a family of random $(n,d)$-graphs. 
\end{thm} 
Theorem~\ref{random_graph} implies that CPs based on random graphs
exhibit fast convergence even when the number of dynamic agents grows unboundedly.
Note that for $n\gg 1$, an $(n,d)$-graph is sparse. Nonetheless, the CP based on 
a random $(n,d)$-graph
possesses the convergence speed that is practically as good 
as that of the normalized complete graph (cf. (\ref{scale_path})).
Therefore, random graphs provide a simple practical way of 
design of CPs that are efficient for coordinating large 
networks. 
\begin{figure}
\begin{center}
{\bf a}\epsfig{figure=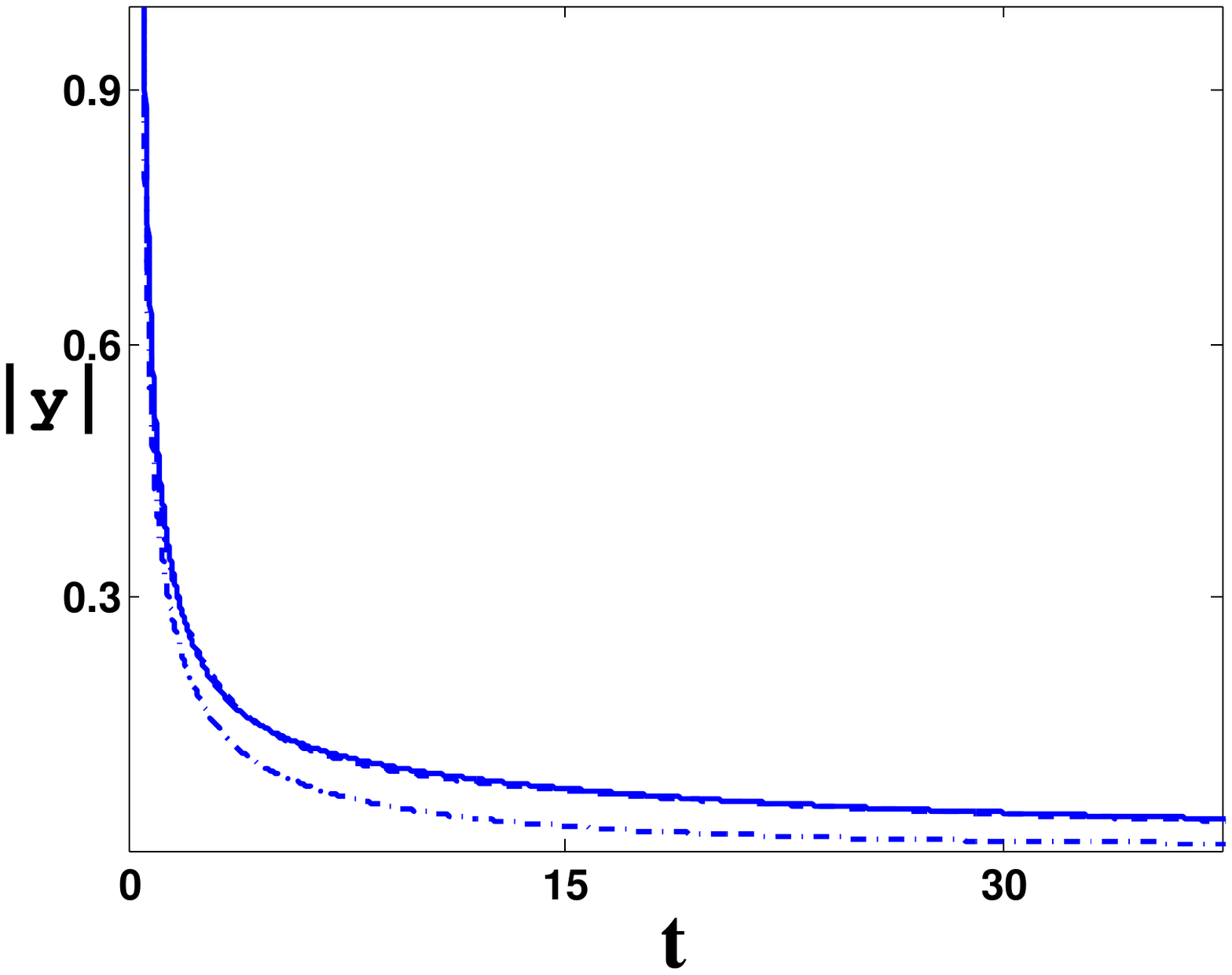, height=2.0in, width=2.5in, angle=0}
\qquad
{\bf b}\epsfig{figure=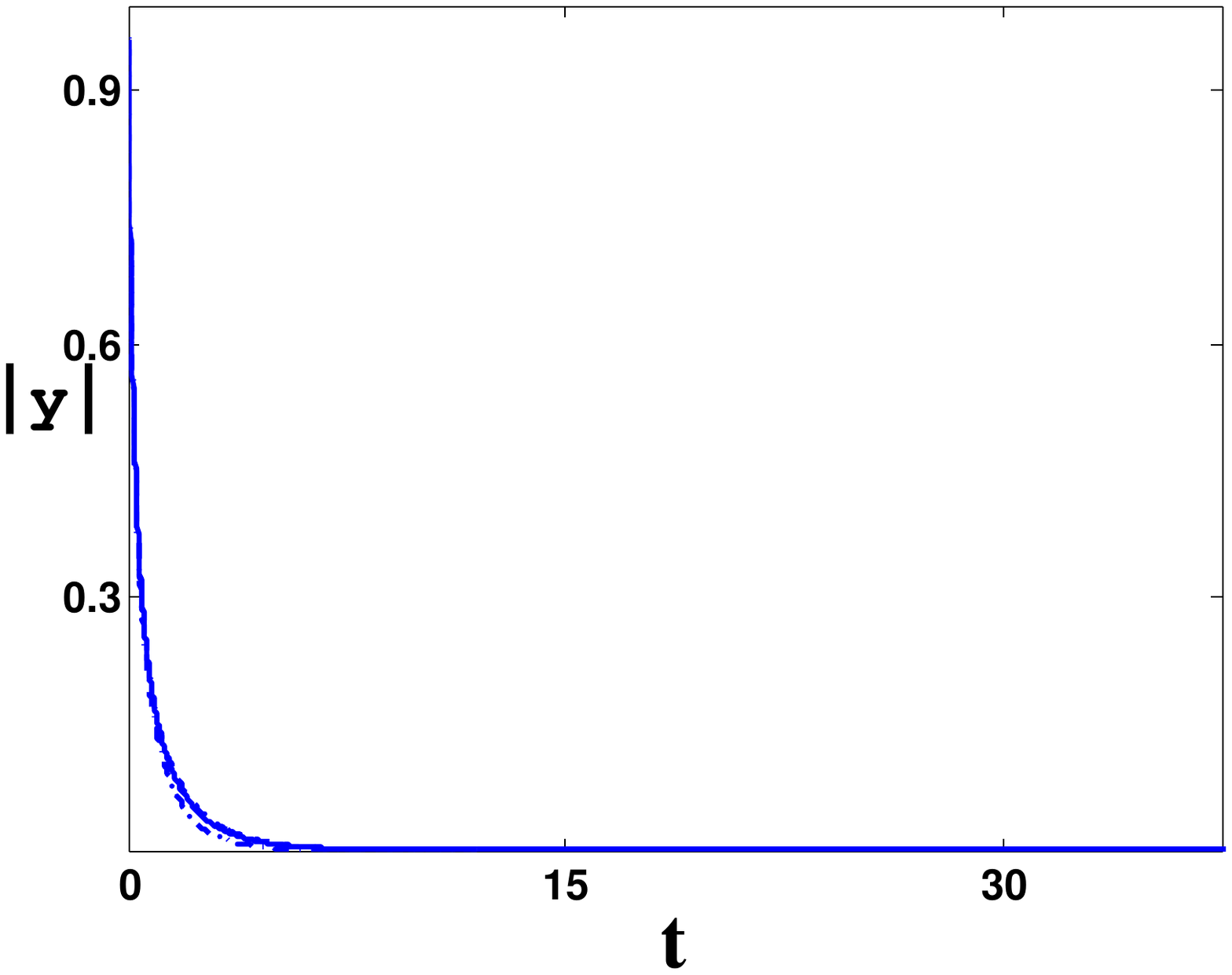, height=2.0in, width=2.5in, angle=0} \\
{\bf c}\epsfig{figure=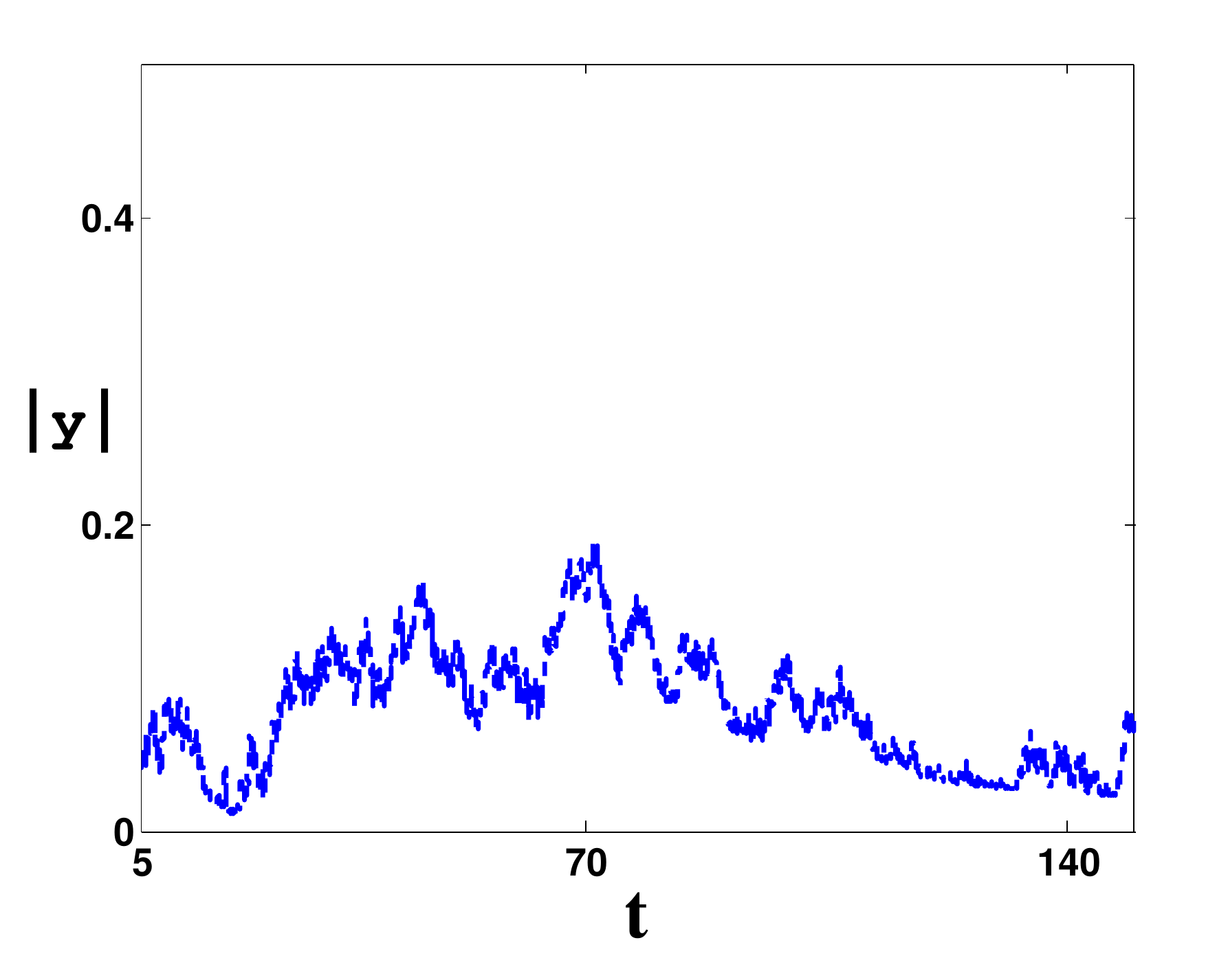, height=2.0in, width=2.5in, angle=0}
\qquad
{\bf d}\epsfig{figure=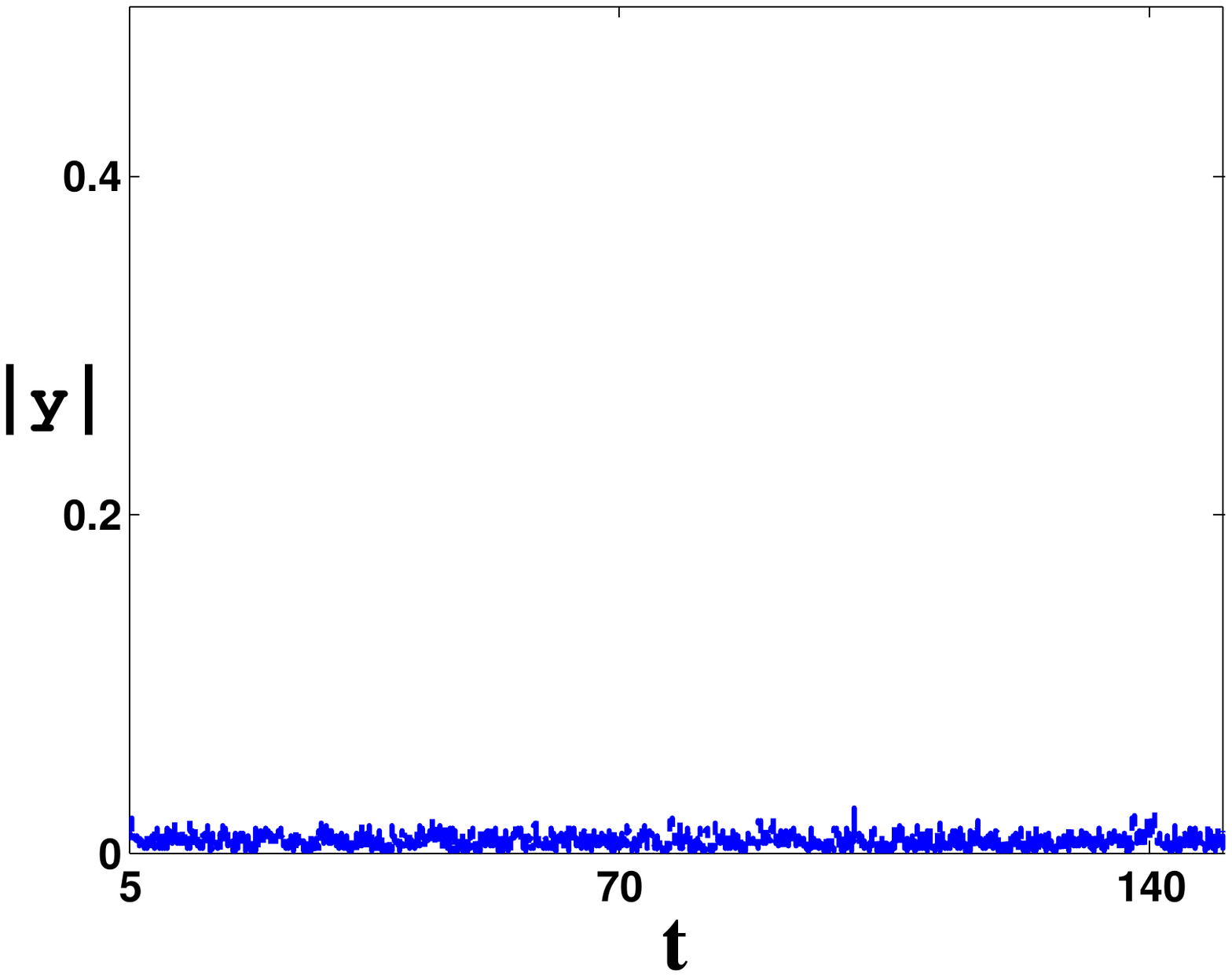, height=2.0in, width=2.5in, angle=0} 
\end{center}
\caption{Numerical simulations of the CPs based on degree $4$
graphs with regular and random connectivity (see Example~\ref{path_and_bipart},
for the definitions of the graphs).
Plots (a) and (b) show convergence of CPs for two networks. For each network
architecture three network sizes were used $100$ (shown in dashed line),
$200$ (dash-dotted line), and $400$ (solid line). The graphs show the
Euclidean norm of the trajectory of the reduced system $y$. The CP based
on random graph in (b) shows much better convergence rate compared to that 
based on a regular graph in (a). Plots in (c) and (d) show the corresponding
results for randomly perturbed CPs. 
}\label{f.2}
\end{figure}

\begin{ex} \lbl{path_and_bipart}
In this example, we compare performance of two CPs 
with regular and random connectivity.
\begin{description}
\item[(a)] The former is a cycle on $n$ vertices, $C_n$.
Each vertex of $C_n$ is connected to $d/2$ ($d$ is even) of 
its nearest neighbors from each side (see Fig.~\ref{f.1}a).
\item[(b)] The latter is a bipartite graph on $2m$ vertices, $B_{2m}$.
The edges are generated using the following algorithm:
\begin{enumerate}
\item
Let $p:[m]\to[m]$ be a random permutation.
In our numerical experiments, we used MATLAB function $\mathsf{randperm}$
to generate random permutations.
For $i\in [m]$, add edge $(i, m+p(i))$.
\item
Repeat step 1. $d-1$ times.
\end{enumerate}
\end{description}
\end{ex}

In Fig.~\ref{f.2}, we present numerical simulations for CPs
based on graphs in Example~\ref{path_and_bipart} for $d=4$
and $n\in\{100, 200, 400\}$. The rates of convergence for
these CPs are summarized in the following table.
\begin{center}
\textsc{
Table 5.4
}
\end{center}
\begin{center}
\begin{tabular}{ |c|| c| c| c| c |}
\hline
$\alpha\backslash n $                & $100$ & $200$ & $400$ & $\infty$\\
\hline
$\alpha(C_n)$ &$.020$ & $.005$ & $.001$ &$0$ \\
$\alpha(B_n)$ & $.597$  &$.554$ & $.547$ & $4-2\sqrt{3}\approx 0.536$\\
\hline
\end{tabular}
\end{center}

The CP based on regular graphs $C_n$ has a very small rate of convergence
already for $n=100$. As $n\to\infty$ $\alpha(C_n)$ tends to zero.
In contrast, random graphs yield rates of convergence with very mild dependence
on the size of the network. For the values of $n$ used in this experiment,
$\alpha(B_n)$ ($n=\{1,2,4\}\times 10^2$) are close to the optimal limiting rate 
$\alpha(B_\infty)=4-2\sqrt{3}\approx 0.536$. The difference of the convergence
rates is clearly seen in the numerical simulations of the corresponding CPs 
(see Fig.~\ref{f.2} a,b). The trajectories generated by CPs with random connections
converge to the consensus subspace faster. We also conducted numerical experiments
with randomly perturbed equation (\ref{perturbed}) to compare the robustness to noise
of the random and regular CPs. The CP on the random graph is more stable to random
perturbations than the one on the regular graph (see Fig.~\ref{f.2}c,d).

\section{Conclusions}\lbl{conclude}
\setcounter{equation}{0}
In this paper, we presented a unified approach to studying convergence
and stochastic stability of a large class of CPs including
CPs on weighted directed graphs; CPs with both positive and negative
conductances, time-dependent CPs, and those under stochastic forcing.
We derived analytical estimates 
characterizing  convergence of CPs and their stability to
random perturbations. Our analysis shows how spectral and 
structural  properties of the graph of the network contribute to 
stability of the corresponding CP.  
In particular, it suggests that the geometry of the cycle subspace
associated with the graph of the CP plays an important role in 
shaping its stability. Further,
we highlighted the advantages of using expanders and,
in particular, random graphs, for CP design. 

The results of this paper elucidate
the link between the structural properties of the graphs and dynamical
performance of CPs. The theory of CPs is closely related to the theory
of synchronization \cite{Blekhman,mmp,pr01,Strogatz03}. With minimal
modifications, the results of the present study carry over to 
a variety coupled one-dimensional dynamical systems 
ranging from the models of power networks \cite{DB10}
to neuronal networks \cite{medvedev09, MZ11, MZ11a}
and drift-diffusion models of decision making \cite{Leonard10a}. 
Moreover, the method of this paper naturally fits in into a general
scheme of analysis of synchronization in coupled systems of multidimensional 
nonlinear oscillators. An interested reader is referred to 
\cite{medvedev10a, medvedev10, MZ11, MZ11a} for related techniques and applications.

\vskip 0.2cm
\noindent
{\bf Acknowledgements.} The author thanks 
Dmitry Kaliuzhnyi-Verbovetskyi for useful discussions
and the referees for careful reading of the manuscript 
and many useful suggestions.
Anatolii Grinshpan provided comments on an earlier version of this
paper, including a more concise proof of Lemma~\ref{spectra}.
Part of this work was done during sabbatical leave at Program of 
Applied and Computational Mathematics of Princeton University. 
This work was supported in part by the NSF Award DMS 1109367.

\end{document}